\documentclass[10pt]{article}
\usepackage[margin=0.67in,top=0.37in,bottom=0.57in]{geometry}
\usepackage{amssymb}
\usepackage{amsmath}
\allowdisplaybreaks
\usepackage{bbm}
\usepackage{enumerate}
\usepackage{float}
\newcommand{\mat}[4]{\left[ \begin{array}{cc} #1 & #2 \\ #3 & #4 \end{array} \right] }
\DeclareMathOperator{\ccc}{circ}
\DeclareMathOperator{\fold}{\bf{Fold}}

\DeclareMathOperator{\unfold}{\bf{Unfold}}

\DeclareMathOperator{\uncirc}{uncirc}
\DeclareMathOperator{\fft}{fft}

\DeclareMathOperator{\ifft}{ifft}

\newcommand{\mary}[2]{{#1}^{\parallel #2}}

\DeclareMathOperator{\rk}{rank}

\DeclareMathOperator{\ind}{ind}
\DeclareMathOperator{\DFT}{\bf{DFT}}
\DeclareMathOperator{\IDFT}{\bf{IDFT}}

\newcommand{\CC}{\mathbbmss{C}}

\usepackage[english, algosection, algoruled, noline]{algorithm2e}

\usepackage{eucal}
\newcommand{\ac}{\mathcal{A}}
\newcommand{\bc}{\mathcal{B}}
\newcommand{\cc}{\mathcal{C}}
\newcommand{\ec}{\mathcal{E}}

\newcommand{\dc}{\mathcal{D}}

\newcommand{\xc}{\mathcal{X}}
\newcommand{\yc}{\mathcal{Y}}
\newcommand{\zc}{\mathcal{Z}}

\newcounter{nex}
\newenvironment{myexample}
  {
   {\noindent \bf Example \thesection.\thenex.}}
   {\stepcounter{nex}
  }

\newtheorem{theorem}{Theorem}[section]
\newtheorem{lemma}{Lemma}[section]
\newtheorem{corollary}{Corollary}[section]

\newtheorem{definition}{Definition}[section]

\numberwithin{equation}{section}
 \date{}

\begin{document}

\title{\bf The generalized inverses of the quaternion tensor via the T-product}

\author{
Hongwei Jin\thanks{School of Mathematics and physics, Guangxi Minzu University, 530006, Nanning, PR China. {\it E-mail address}:  jhw-math@126.com} \quad
Peifeng Zhou\thanks{School of Mathematics and physics, Guangxi Minzu University, 530006,
Nanning, PR China. {\it E-mail address}:  peifengzhou1004@126.com}\quad
Hongjie Jiang\thanks{Corresponding author. School of Mathematics and physics, Guangxi Minzu University, 530006,
Nanning, PR China. {\it E-mail address}:  hongjiejiang@yeah.net}\quad
Xiaoji Liu\thanks{School of Mathematics and physics, Guangxi Minzu University, 530006,
Nanning, PR China. {\it E-mail address}:  xiaojiliu72@126.com}
\quad}

\maketitle

\begin{abstract}
 In this article, specific definitions of the Moore-Penrose inverse, Drazin inverse of the quaternion tensor and the inverse along two quaternion tensors are introduced under the T-product. Some characterizations, representations and properties of the defined inverses are investigated. Moreover, algorithms are established for computing the Moore-Penrose inverse, Drazin inverse of the quaternion tensor and the inverse along two quaternion tensors, respectively.\\



\end{abstract}

\section{Introduction}
We denote by $\mathbb{H}$ the algebra of the quaternion, which is introduced by Hamilton in 1843.
An element $q$ of $\mathbb{H}$ is of the form
\begin{equation*}
  q = a + b{\bf{i}} + c{\bf{j}} + d{\bf{k}}, a, b, c, d\in\mathbb{R},
\end{equation*}
where ${\bf{i}}, {\bf{j}}$ and ${\bf{k}}$ are imaginary units. By definition, they satisfy
$${\bf{i}}^2 = {\bf{j}}^2 = {\bf{k}}^2 = {\bf{ijk}} = -1.$$
Given $q = a + b{\bf{i}} + c{\bf{j}} + d{\bf{k}}$, then
\begin{itemize}
\item [] (i)  the conjugate quaternion of $q$ is $\overline{q}=a - b{\bf{i}} - c{\bf{j}} - d{\bf{k}}$

\item [] (ii) the norm of $q$ is $|q|=\sqrt{q\overline{q}}=\sqrt{a^2+b^2+c^2+d^2}$

\item  [] (iii) the real and the imaginary parts of $q$ are respectively $\text{Re}(q) = \frac{1}{2}(q+\overline{q}) = a$ and $\text{Im}(q) = \frac{1}{2} (q-\overline{q}) = b{\bf{i}} + c{\bf{j}} + d{\bf{k}}$.
\end{itemize}

It is well know that the quaternion is not only part of
contemporary mathematics \cite{Zhang1,Huang1,Zhang2,FL1}, but also widely and heavily used in computer graphics,
control theory, quantum physics, signal and color image processing, and so on \cite{LS1,BS1,SB1}.

An order $3$ tensor $\mathcal{A}=(a_{i_1}a_{i_2}a_{i_3} )$, $1\leq i_j\leq n_j$, $(j = 1, 2, 3)$ is a multidimensional array with
$n_1n_2n_3$ entries. Let $\mathbb{R}^{n_1\times n_2\times n_3}$ , $\mathbb{C}^{n_1\times n_2\times n_3}$ and $\mathbb{H}^{n_1\times n_2\times n_3}$ stand, respectively, for the sets of
the order 3 tensors over the real number field $\mathbb{R}$, the complex number
field $\mathbb{C}$ and the real quaternion algebra
$$\mathbb{H} = \{a + b\textbf{i} + c\textbf{j} + d\textbf{k} \mid \textbf{i}^2 = \textbf{j}^2 = \textbf{k}^2 = \textbf{ijk} = -1, a, b, c, d\in \mathbb{R}\}.$$

Higher-order tensors arise in a wide variety of
application areas, including psychometrics \cite{Kroonenberg}, chemometrics \cite{Smilde}, image
and signal processing \cite{Comon,Lathauwer,Nagy,Sidiropoulos,Hoge,Rezghi,HKBH} and so on.

The generalized inverse of an arbitrary matrix
has many applications in statistics, prediction theory, control system analysis, curve fitting,
numerical analysis and the solution of linear integral equations \cite{B-G,CamMey,WWQ,B,Rao6,Hunter,Puri,Tran}.

The generalized inverses of tensors based on different tensor products have
been investigated frequently. Sun et al.\cite{Sun2} defined the Moore-Penrose inverse of tensors with the Einstein
product, and the explicit formulas of the Moore-Penrose inverse of some block
tensors were obtained. Krushnachandra et al.\cite{Krushnachandra} proved some more identities involving the Moore-Penrose inverse of tensors and   obtained a few necessary and sufficient conditions for the reverse order law for the Moore-Penrose inverse of tensors via the Einstein product.   Sahoo et al.\cite{Sahoo} introduced  specific definitions of the core and core-EP inverses of complex tensors. Some characterizations, representations and properties of the core and core-EP inverses were investigated.
Miao et al.\cite{Miao}  presented the definition of generalized tensor function according to the tensor singular value decomposition (T-SVD) via the tensor T-product.
Miao et al.\cite{Miao2}  proposed the T-Jordan canonical form and  studied the T-group inverse
and T-Drazin inverse.
Ji  et al.\cite{Ji} defined the null spaces and the ranges of tensors, and study their relationship. The fundamental theorem of linear algebra for matrix spaces to tensor spaces were extended.
Ji  et al.\cite{Ji2} had a further study the properties of even-order tensors with Einstein product. The notion of the Drazin inverse of a square matrix to an even-order square tensor were extended. An expression for the Drazin inverse through the core-nilpotent decomposition for a tensor of even-order was obtained. Jin et al.\cite{Jin}  established the generalized inverse
of tensors by using tensor equations. Moreover, the authors investigated the least squares
solutions of tensor equations. He et al.\cite{He} investigated and discussed in detail the structures of
quaternion tensor SVD, quaternion tensor rank decomposition, and $\eta$-Hermitian
quaternion tensor decomposition with the isomorphic group structures and Einstein product. Stanimirovic et al.\cite{S} considered computation of tensor outer with prescribed range and kernel
of higher order tensors. Further,  conditions
for the existence, representation and computation of the Moore-Penrose inverse, the weighted Moore-Penrose inverse, the Drazin
inverse and the usual inverse of tensors are derived.

We organize the paper as follows. In the next subsection, we introduce some notations and
definitions which are helpful in proving the main results. In Section 3, we
will firstly define the Moore-Penrose inverse of the quaternion tensor. Then, an expression of the Moore-Penrose inverse is given by using the quaternion tensor SVD. In Section 4, we will establish an algorithm for computing the Moore-Penrose inverse of the quaternion tensor. In Section 5, we will be concerned on the Drazin inverse of the quaternion tensor. We investigate definitions, characterizations, representations and
properties for this inverse.  Section 6 contains the context of the inverse along two quaternion tensors. The definitions of the right and left inverse along two quaternion tensors are given. Some expressions of the right and left inverse along two quaternion tensors are obtained. In the final section, we get an algorithm for computing the right/left inverse along two quaternion tensors.

\section{Preliminaries}

Throughout this paper tensors are denoted by Euler script letters (e.g., $\mathcal{A}$, $\mathcal{B}$, $\mathcal{C}$,...), while capital letters represent matrices, boldface lowercase letters represent vectors, and lowercase letters refer to scalars. $H^n, H^{m\times n}, H_r^{m\times n}$ stand, respectively, for the sets of quaternion vector space, the quaternion matrix space and the quaternion matrix space with rank $r$. $\textbf{a}^H$, $A^H$ and $\mathcal{A}^H$ are the conjugate transpose of $\textbf{a}$, $A$ and $\mathcal{A}$, respectively.

Now, we present the T-product defined in the work \cite{KM,MSL}.





















Let $\textbf{a}\in{{H}}^n$. Recall that if $\textbf{a}=\begin{bmatrix}a_1 & a_2 & ... & a_n\end{bmatrix}^H$, then
\begin{equation*}
  \ccc(\textbf{a})=\begin{bmatrix}
                     a_1 & a_n & \cdots & a_2 \\
                     a_2 & a_1 & \cdots & a_3 \\
                     \vdots & \vdots &  & \vdots \\
                     a_n & a_{n-1} & \cdots & a_1 \\
                   \end{bmatrix}
\end{equation*}
is a circulant quaternion matrix. Similarly, if $A_1,...,A_n$ are $n_1\times n_2$ quaternion matrices, then
\begin{equation*}
  \ccc(A_1,...,A_n)=\begin{bmatrix}
                     A_1 & A_n & \cdots & A_2 \\
                     A_2 & A_1 & \cdots & A_3 \\
                     \vdots & \vdots &  & \vdots \\
                     A_n & A_{n-1} & \cdots & A_1 \\
                   \end{bmatrix}.
\end{equation*}

Define $\unfold(\cdot)$ to take an $n_1\times n_2\times n_3$ quaternion tensor and return an $n_1n_3\times n_2$ block quaternion matrix in the following way:
\begin{equation*}\label{uu}
\unfold(\mathcal{A})=\begin{bmatrix}
\overline{{A}}_1 \\
\overline{{A}}_2 \\
\vdots \\
\overline{{A}}_{n_3} \\
\end{bmatrix},
\end{equation*}
where $\overline{{A}}_1=\mathcal{A}(:,:,1)$,\ldots,$\overline{{A}}_{n_3}=\mathcal{A}(:,:,n_3)$, $\mathcal{A}(:,:,i)$, $i=1,2,...,n_3$ are the frontal slices of $\mathcal{A}$. $\fold(\cdot)$ is the inverse operation, which takes
an $n_1n_3\times n_2$ block quaternion matrix and returns an $n_1\times n_2\times n_3$ quaternion tensor. In additon, we have
\begin{equation*}\label{aa15}
  \fold(\unfold(\mathcal{A}))=\mathcal{A}.
\end{equation*}

Now, one can create a quaternion tensor in a block circulant pattern, where each block is a quaternion matrix:
\begin{equation}\label{zzx2}
  \ccc(\unfold(\mathcal{A}))=\begin{bmatrix}
                      \overline{{A}}_1 & \overline{{A}}_{n_3} & \overline{{A}}_{n_3-1} & \cdots & \overline{{A}}_2  \\
                      \overline{{A}}_2 & \overline{{A}}_1 & \overline{{A}}_{n_3} & \cdots & \overline{{A}}_3\\
                      \vdots & \vdots & \vdots &  & \vdots \\
                      \overline{{A}}_{n_3} & \overline{{A}}_{n_3-1} & \overline{{A}}_{n_3-2} & \cdots &\overline{{A}}_1 \\
                    \end{bmatrix}.
\end{equation}


Next, we will use the Kronecker product, symbolized as $\otimes$. Let $F_{n_i}$ be the $n_i\times n_i$
discrete Fourier transform (DFT) matrix, i.e.,
\begin{equation*}
  F_{n_i}=\frac{1}{\sqrt{n_i}}\left[
                                \begin{array}{cccccc}
                                  1 & 1 & 1 & 1 & \cdots & 1 \\
                                  1 & \omega & \omega^2 & \omega^3 & \cdots & \omega^{(n_i-1)} \\
                                  1 & \omega^2 & \omega^4 & \omega^6 & \cdots & \omega^{2(n_i-1)} \\
                                  1 & \omega^3 & \omega^6 & \omega^9 & \cdots & \omega^{3(n_i-1)} \\
                                  \vdots & \vdots & \vdots & \vdots &  & \vdots \\
                                  1 & \omega^{(n_i-1)} & \omega^{2(n_i-1)} & \omega^{3(n_i-1)} & \cdots & \omega^{(n_i-1)(n_i-1)} \\
                                \end{array}
                              \right],
\end{equation*}
where $\omega=e^{2\pi i/n_i}$. Define $F = F_{n_3}$. Then there exist quaternion matrices $A_1, \ldots, A_{n_3}
\in{H}^{n_1 \times n_2}$ such that

\begin{equation} \label{Faa}
(F \otimes I_{n_1}) \ccc(\unfold(\mathcal{A})) (F^H\otimes I_{n_2}) =
\begin{bmatrix} A_1 & 0 & \cdots & 0 \\ 0 & A_2 & \cdots & 0 \\
\vdots & \vdots & \ddots & \vdots  \\ 0 & 0 & \cdots   &  A_{n_3}
\end{bmatrix}.
\end{equation}

Let us define the function $\DFT(\cdot)$ that implements the above procedure, i.e,
\begin{equation*}
\textbf{DFT}(\ccc(\unfold(\mathcal{A})))=(F\otimes I_{n_1})\ccc(\unfold(\mathcal{A}))(F^H\otimes I_{n_2}).
\end{equation*}
Also, let $\text{uncirc}(\cdot)$ and $\textbf{IDFT}(\cdot)$ be the inverse operations of $\text{circ}(\cdot)$ and $\textbf{DFT}(\cdot)$, respectively.

The formula (\ref{zzx2}) allows us to define the T-product of two quaternion tensors.

\begin{definition}{\rm\cite{KM,MSL}}\label{d1}
Let $\mathcal{A}\in\mathbb{H}^{n_1\times n_2\times n_3}$ and $\mathcal{B}\in\mathbb{H}^{n_2\times l\times n_3}$.
Then the {\bf{T-product}} $\mathcal{A}*\mathcal{B}$ is the $n_1\times l\times n_3$ quaternion tensor defined as
\begin{equation}\label{dd1}
  \mathcal{A}*\mathcal{B}=\fold(\ccc(\unfold(\mathcal{A}))\unfold(\mathcal{B})).
\end{equation}
\end{definition}

How to compute this new product with mathematical software?  In the following, an iterative procedure
in MATLAB pseudocode is provided which is inspired by \cite{MSL}.

\begin{algorithm}[H]
\caption{\textsc{Compute the T-product of two quaternion tensors }}\label{algo:sylvester9}
\KwIn{$n_1\times n_2\times n_3$ quaternion tensor $\mathcal{A}$ and $n_2\times l\times n_3$ quaternion tensor $\mathcal{B}$}
\KwOut{$n_1\times l\times n_3$ quaternion tensor $\mathcal{C}$}
\begin{enumerate}
\addtolength{\itemsep}{-0.8\parsep minus 0.8\parsep}

\item $i=3$, $\widetilde{\mathcal{A}} = \fft(\mathcal{A},[ ~ ],i)$, $\widetilde{{\mathcal{B}}} = \fft(\mathcal{B},[ ~ ],i)$

\item  for $i=1,\ldots, n_3$

\quad $\widetilde{{\mathcal{C}}}(:,:,i)=
\widetilde{{\mathcal{A}}}(:,:,i)\cdot\widetilde{{\mathcal{B}}}(:,:,i)$,

end

\item  $i=3$, $\mathcal{C}=\ifft(\widetilde{{\mathcal{C}}},[ ~ ],i)$.
\end{enumerate}
\end{algorithm}
Denote
\begin{equation*}
\textbf{DFT}(\ccc(\unfold(\mathcal{A})))=
\begin{bmatrix} A_1 & 0 & \cdots & 0 \\ 0 & A_2 & \cdots & 0 \\
\vdots & \vdots & \ddots & \vdots  \\ 0 & 0 & \cdots   &  A_{n_3}
\end{bmatrix}=\widehat{A}.
\end{equation*}

By Algorithm \ref{algo:sylvester9}, we have the following fact:
\begin{equation}\label{a9}
  \mathcal{A}*\mathcal{B}=\mathcal{C}\Longleftrightarrow\widehat{A}\cdot\widehat{B}=\widehat{C},
\end{equation}
where $\widehat{B}$ and $\widehat{C}$ can be obtained by using the same manner of $\widehat{A}$.

Now, it is easy to check the following basic properties of the T-product.
\begin{lemma}\label{ppa}
If $\mathcal{A}, \mathcal{B}, \mathcal{C}$ are quaternion tensors of adequate size, then the following statements are true:
\begin{itemize}
\item [] $(i)$   $\mathcal{A}*(\mathcal{B}+\mathcal{C})=\mathcal{A}*\mathcal{B}+\mathcal{A}*\mathcal{C}$;
\item [] $(ii)$  $(\mathcal{A}+\mathcal{B})*\mathcal{C}=\mathcal{A}*\mathcal{C}+\mathcal{B}*\mathcal{C}$;
\item [] $(iii)$   $(\mathcal{A}*\mathcal{B})*\mathcal{C}=\mathcal{A}*(\mathcal{B}*\mathcal{C})$.
\end{itemize}
\end{lemma}


\begin{definition}{\rm\cite{MSL}}\label{d5}
The $n\times n\times n_3$  {\bf {identity tensor}} $\mathcal{I}$ is
the tensor such that the first frontal
slice $\mathcal{I}(:,:,1)$ is the $n\times n$ identity matrix and all other frontal slices $\mathcal{I}(:,:,k)$, $k=2,...,n_3$ are zero matrices.
\end{definition}

The following result can be checked by using (\ref{a9}).

\begin{lemma}\label{pp2a}
Let $\mathcal{A}\in\mathbb{H}^{n\times n\times n_3}$ and $\mathcal{I}$ be an $n\times n\times n_3$  order-$3$ identity tensor. Then, $$\mathcal{I}*\mathcal{A}=\mathcal{A}*\mathcal{I}=\mathcal{A}.$$
\end{lemma}
{\sc Proof:}
Let
\begin{equation*}
\textbf{DFT}(\ccc(\unfold(\mathcal{I})))=
\begin{bmatrix} I_1 & 0 & \cdots & 0 \\ 0 & I_2 & \cdots & 0 \\
\vdots & \vdots & \ddots & \vdots  \\ 0 & 0 & \cdots   &  I_{n_3}
\end{bmatrix}=\widehat{I} \ \text{and} \ \textbf{DFT}(\ccc(\unfold(\mathcal{A})))=
\begin{bmatrix} A_1 & 0 & \cdots & 0 \\ 0 & A_2 & \cdots & 0 \\
\vdots & \vdots & \ddots & \vdots  \\ 0 & 0 & \cdots   &  A_{n_3}
\end{bmatrix}=\widehat{A}.
\end{equation*}
Since $\widehat{I}\cdot\widehat{{A}}=\widehat{A}$, using (\ref{a9}), we have
\begin{eqnarray*}
  \widehat{I}\cdot\widehat{{A}}=\widehat{A}&\Rightarrow&\mathcal{I}*\mathcal{A}=\mathcal{A}.
\end{eqnarray*}
Similarly, $\mathcal{A}*\mathcal{I}=\mathcal{A}$. $\Box$

\begin{definition}\label{d9}
Let $\mathcal{A}\in\mathbb{H}^{n\times n\times n_3}$. If there exists an tensor $\mathcal{B}\in\mathbb{H}^{n\times n\times n_3}$ such that
\begin{equation*}\label{0}
  \mathcal{A}*\mathcal{B}=\mathcal{I} \ \  \text{and}    \  \ \mathcal{B}*\mathcal{A}=\mathcal{I},
\end{equation*}
then $\mathcal{A}$ is said to be {\bf{invertible}}. Moreover, $\mathcal{B}$ is the {\bf{inverse}} of $\mathcal{A}$, which is denoted by $\mathcal{A}^{-1}$.
\end{definition}

In fact, the inverse of an invertible quaternion tensor is unique.

\begin{lemma}\label{pp5a}
If $\mathcal{A}\in\mathbb{H}^{n\times n\times n_3}$ is invertible, then its inverse tensor is unique.
\end{lemma}
{\sc Proof:}
Suppose $\mathcal{B}_1$ and $\mathcal{B}_2$ are both the inverse of $\mathcal{A}$. Since $\mathcal{A}$ is invertible, one has
\begin{equation*}
  \mathcal{B}_1=\mathcal{B}_1*\mathcal{I}=\mathcal{B}_1*(\mathcal{A}*\mathcal{B}_2)=(\mathcal{B}_1*\mathcal{A})*\mathcal{B}_2=\mathcal{I}*\mathcal{B}_2=\mathcal{B}_2,
\end{equation*}
where we used Lemma \ref{ppa} (c), Lemma \ref{pp2a} and Definition \ref{d9}.
Hence, the inverse tensor of $\mathcal{A}$ is unique.
$\Box$

The conjugate transpose of quaternion tensors can be defined as follows.

\begin{definition}\label{d92}
If $\mathcal{A}\in\mathbb{H}^{n_1\times n_2\times n_3}$, then the {\bf{conjugate  transpose}} of $\mathcal{A}$, which is denoted by $\mathcal{A}^H$, is the $n_2\times n_1\times n_3$ quaternion tensor obtained by conjugate transposing each $\mathcal{A}(:,:,i)$, for $i=1,2,\ldots,n_3$ and then reversing the order of the $\mathcal{A}(:,:,i)$ from $2$ through $n_3$, i.e.,
\begin{equation*}
 \mathcal{A}^H(:,:,1)=(\mathcal{A}(:,:,1))^H,
\end{equation*}
\begin{equation}\label{a01}
 \mathcal{A}^H(:,:,i)=(\mathcal{A}(:,:,n_3+2-i))^H, \ i=2,\ldots,n_3.
\end{equation}
\end{definition}

\begin{lemma}\label{ll}
Suppose that $\mathcal{A}$, $\mathcal{B}$ are two quaternion tensors such that $\mathcal{A}*\mathcal{B}$ and $\mathcal{B}^H*\mathcal{A}^H$ are
defined. Then,
\begin{itemize}
\item [] $(i)$ $(\mathcal{A}^H)^H=\mathcal{A}$.
\item [] $(ii)$ $(\mathcal{A}+\mathcal{B})^H=\mathcal{A}^H+\mathcal{B}^H$.
\item [] $(iii)$ $(\mathcal{A}*\mathcal{B})^H=\mathcal{B}^H*\mathcal{A}^H$.
\end{itemize}
\end{lemma}
{\sc Proof:}
(i) By (\ref{a01}), we have
$$(\mathcal{A}^H)^H(:,:,1)=((\mathcal{A}(:,:,1))^H)^H=\mathcal{A}(:,:,1).$$
$$(\mathcal{A}^H)^H(:,:,i)=((\mathcal{A}(:,:,i))^H)^H=\mathcal{A}(:,:,i),\ i=2,\ldots,n_3.$$
Hence, $(\mathcal{A}^H)^H=\mathcal{A}$.

(ii) Again using (\ref{a01}), one has
\begin{equation*}
(\mathcal{A}+\mathcal{B})^H(:,:,1)=(\mathcal{A}(:,:,1)+\mathcal{B}(:,:,1))^H
=(\mathcal{A}(:,:,1))^H+(\mathcal{B}(:,:,1))^H=\mathcal{A}^H(:,:,1)+\mathcal{B}^H(:,:,1).
\end{equation*}
\begin{eqnarray*}
  (\mathcal{A}+\mathcal{B})^H(:,:,i)&=&(\mathcal{A}(:,:,n_3+2-i)+\mathcal{B}(:,:,n_3+2-i))^H\\
&=&(\mathcal{A}(:,:,n_3+2-i))^H+(\mathcal{B}(:,:,n_3+2-i))^H\\
 &=& \mathcal{A}^H(:,:,n_3+2-i)+\mathcal{B}^H(:,:,n_3+2-i),
\end{eqnarray*}
where $i=2,\ldots,n_3.$ Hence, $(\mathcal{A}+\mathcal{B})^H=\mathcal{A}^H+\mathcal{B}^H$.

(iii) Suppose $(\mathcal{A}*\mathcal{B})^H=\mathcal{C}^H$, by (\ref{a9}), we have
\begin{equation*}
  (\mathcal{A}*\mathcal{B})^H=\mathcal{C}^H\Rightarrow(\widehat{A}\cdot\widehat{B})^H=\widehat{C}^H
  \Rightarrow\widehat{C}^H=\widehat{B}^H\cdot\widehat{A}^H\Rightarrow\mathcal{C}^H
  =\mathcal{B}^H*\mathcal{A}^H.
\end{equation*}
Hence, $(\mathcal{A}*\mathcal{B})^H=\mathcal{B}^H*\mathcal{A}^H$. $\Box$


The following definitions are useful in establishing the main results.



\begin{definition}\label{d10}
An $n\times n\times n_3$ quaternion tensor $\mathcal{Q}$ is
{\bf{orthogonal}} if
$$\mathcal{Q}^H*\mathcal{Q}=\mathcal{Q}*\mathcal{Q}^H=\mathcal{I}.$$
\end{definition}



\begin{definition}
Let $\mathcal{A}=(a_{i_1i_2i_3})\in\mathbb{H}^{n_1\times n_2\times n_3}$. Then, $\mathcal{A}$ is called an {\bf{F-diagonal}} quaternion tensor if $a_{i_1i_2i_3}=0$ when $i_1\neq i_2$.
\end{definition}

\section{The Moore-Penrose Inverse of Quaternion Tensors}

The T-product of two tensors presented in Definition \ref{d1} allows us to obtain the Moore-Penrose inverse of an arbitrary quaternion tensor $\mathcal{A}$.

\begin{definition}\label{d4}
Let $\mathcal{A}\in\mathbb{H}^{n_1\times n_2\times n_3}$. If there exists a quaternion tensor $\mathcal{X}\in\mathbb{H}^{n_2\times n_1\times n_3}$ such that
\begin{equation}\label{mp}
(1) \, \, \mathcal{A}*\mathcal{X}*\mathcal{A}=\mathcal{A} \qquad
(2) \, \, \mathcal{X}*\mathcal{A}*\mathcal{X}=\mathcal{X} \qquad
(3) \, \, (\mathcal{A}*\mathcal{X})^{H}=\mathcal{A}*\mathcal{X} \qquad
(4) \, \, (\mathcal{X}*\mathcal{A})^{H}=\mathcal{X}*\mathcal{A},
\end{equation}
then $\mathcal{X}$ is called the {\bf{Moore-Penrose inverse}} of the quaternion tensor $\mathcal{A}$ and is denoted by $\mathcal{A}^\dag$.
\end{definition}

For any $\mathcal{A}\in\mathbb{H}^{n_1\times n_2\times n_3}$,  denote $\mathcal{A}{\{i, j,  \ldots,  k\}}$
the set of all  $\mathcal{X}\in\mathbb{H}^{n_2\times n_1\times n_3}$ which satisfy equations ($i$),  ($j$), $\ldots$ , ($k$) of $(\ref{mp})$.  In this case, $\mathcal{X}$ is a $\{i, j, \ldots, k\}$-inverse.

If $\mathcal{A}$ is invertible, it is clear that $\mathcal{X}=\mathcal{A}^{-1}$ trivially satisfies the four equations.

\begin{theorem}{\rm\cite{Zhang1}} (Quaternion matrix SVD) Let $A\in{H^{m\times n}}$
be of rank $r$. Then there exist
unitary quaternion matrices $U\in{H^{m\times m}}$ and $V\in{H^{n\times n}}$
such that
\begin{equation*}
  A=U\begin{bmatrix}
  \Sigma_r & 0   \\
  0 & 0   \\
\end{bmatrix}V^H,
\end{equation*}
where $\Sigma_r=diag(\sigma_1,\ldots,\sigma_r)$ and the $\sigma_i$, $i=1,\ldots,r$ are real positive singular values of $A$.
\end{theorem}

\begin{theorem}\label{vt}(Quaternion tensor SVD) Let $\mathcal{A}\in\mathbb{H}^{n_1\times n_2\times n_3}$. Then there exist
unitary quaternion tensor $\mathcal{U}\in\mathbb{H}^{n_1\times n_1\times n_3}$ and $\mathcal{V}\in\mathbb{H}^{n_2\times n_2\times n_3}$
such that
\begin{equation*}
  \mathcal{A}=\mathcal{U}*\mathcal{S}*\mathcal{V}^H,
\end{equation*}
where $\mathcal{S}$ is an $n_1\times n_2\times n_3$ {$F$-diagonal}  tensor.
\end{theorem}
{\sc Proof:} Let
\begin{eqnarray*} 
\textbf{DFT}(\text{circ}(\textbf{Unfold}(\mathcal{A}))) =
\begin{bmatrix} A_1 & & \\ & \ddots & \\ & & A_{n_3} \end{bmatrix}.
\end{eqnarray*}
Next compute the SVD of each $A_i$, $i=1,2,\ldots,n_3$ as $A_i = U_i\Sigma_iV^H_i$. Thus, we have
\begin{equation}\label{oooo}
  \begin{bmatrix}
                     A_1 &   &   \\
                       & \ddots   &   \\
                       &   &  A_{n_3} \\
                   \end{bmatrix}=\begin{bmatrix}
                     U_1 &   &   \\
                       & \ddots   &   \\
                       &   &  U_{n_3} \\
                   \end{bmatrix}\begin{bmatrix}
                     \Sigma_1 &   &   \\
                       & \ddots   &   \\
                       &   &  \Sigma_{n_3} \\
                   \end{bmatrix}\begin{bmatrix}
                     V^H_1 &   &   \\
                       & \ddots   &   \\
                       &   &  V^H_{n_3} \\
                   \end{bmatrix}.
\end{equation}
Then,
\begin{eqnarray*}\label{oooo2}
  \mathcal{A}&=&\textbf{Fold}(\text{uncirc}(\textbf{IDFT}(\begin{bmatrix}
                     A_1 &   &   \\
                       & \ddots   &   \\
                       &   &  A_{n_3} \\
                   \end{bmatrix})))\\
  &=&\textbf{Fold}(\text{uncirc}(\textbf{IDFT}(\begin{bmatrix}
                     U_1 &   &   \\
                       & \ddots   &   \\
                       &   &  U_{n_3} \\
                   \end{bmatrix})))*\textbf{Fold}(\text{uncirc}(\textbf{IDFT}(\begin{bmatrix}
                     \Sigma_1 &   &   \\
                       & \ddots   &   \\
                       &   &  \Sigma_{n_3} \\
                   \end{bmatrix})))*\\
                   &&\textbf{Fold}(\text{uncirc}(\textbf{IDFT}(\begin{bmatrix}
                     V^H_1 &   &   \\
                       & \ddots   &   \\
                       &   &  V^H_{n_3} \\
                   \end{bmatrix})))\\
                   &=&\mathcal{U}*\mathcal{S}*\mathcal{V}^H.
\end{eqnarray*}
This completes the proof. $\Box$

In the following, we will show the existence and uniqueness of the Moore-Penrose inverse of a quaternion tensor $\mathcal{A}$.





\begin{theorem}\label{t2}
The Moore-Penrose inverse of an arbitrary quaternion tensor $\mathcal{A}\in\mathbb{H}^{n_1\times n_2\times n_3}$ exists and is unique.
\end{theorem}
{\sc Proof:} Let
\begin{eqnarray*} 
\textbf{DFT}(\text{circ}(\textbf{Unfold}(\mathcal{A}))) =
\begin{bmatrix} A_1 & & \\ & \ddots & \\ & & A_{n_3} \end{bmatrix}
\quad \text{and} \quad
\textbf{DFT}(\text{circ}(\textbf{Unfold}(\mathcal{X}))) =
\begin{bmatrix} X_1 & & \\ & \ddots & \\ &  &  X_{n_3} \end{bmatrix}.
\end{eqnarray*}
Suppose $A_i=U_i\Sigma_iV^H_i$ are the SVD of each $A_i$, $i=1,...,n_3$ and for each $\Sigma_i=diag(\sigma^i_{1},...,\sigma^i_{r_i},0,...,0)$, where $r_i=rank(A_i)$, $\sigma^i_{1},...,\sigma^i_{r_i}$ are real positive singular values of $A_i$.
We define the matrices $R_i=diag(\frac{1}{\sigma^i_{1}},...,\frac{1}{\sigma^i_{r_i}},0,...,0)$, for $i=1,...,{n_3}$.
Observe that $R_i = \Sigma_i^\dag$. Let
$X_i = V_i R_i U_i^H$ for $i=1, \ldots, {n_3}$. Now, we have
\begin{equation}\label{svdd}
  \begin{bmatrix}
                     X_1 &   &   \\
                       & \ddots   &   \\
                       &   &  X_{n_3} \\
                   \end{bmatrix}=\begin{bmatrix}
                     V_1 &   &   \\
                       & \ddots   &   \\
                       &   &  V_{n_3} \\
                   \end{bmatrix}\begin{bmatrix}
                     R_1 &   &   \\
                       & \ddots   &   \\
                       &   &  R_{n_3} \\
                   \end{bmatrix}\begin{bmatrix}
                     U^H_1 &   &   \\
                       & \ddots   &   \\
                       &   &  U^H_{n_3} \\
                   \end{bmatrix}.
\end{equation}
Then,
\begin{eqnarray*}\label{oooo2}
  \mathcal{X}&=&\fold(\text{uncirc}(\textbf{IDFT}(\begin{bmatrix}
                     X_1 &   &   \\
                       & \ddots   &   \\
                       &   &  X_{n_3} \\
                   \end{bmatrix})))\\
  &=&\fold(\text{uncirc}(\textbf{IDFT}(\begin{bmatrix}
                     V_1 &   &   \\
                       & \ddots   &   \\
                       &   &  V_{n_3} \\
                   \end{bmatrix})))*\fold(\text{uncirc}(\textbf{IDFT}(\begin{bmatrix}
                     R_1 &   &   \\
                       & \ddots   &   \\
                       &   &  R_{n_3} \\
                   \end{bmatrix})))*\\
                   &&\fold(\text{uncirc}(\textbf{IDFT}(\begin{bmatrix}
                     U^H_1 &   &   \\
                       & \ddots   &   \\
                       &   &  U^H_{n_3} \\
                   \end{bmatrix})))\\
                   &=&\mathcal{V}*\mathcal{R}*\mathcal{U}^H.
\end{eqnarray*}
where $\mathcal{U}$, $\mathcal{V}$ are orthogonal $n_1\times n_1\times n_3$, $n_2\times n_2\times n_3$ quaternion tensors, respectively, and $\mathcal{R}$ is an $n_2\times n_1\times n_3$ $F$-diagonal tensor.  One can check that $\mathcal{X}$ satisfies $(\ref{mp})$.

On the other hand, let $\mathcal{X}_1$ and $\mathcal{X}_2$ be  solutions of $(\ref{mp})$. One has
\begin{eqnarray*}
  \mathcal{X}_1 &=& \mathcal{X}_1*\mathcal{A}*\mathcal{X}_1=\mathcal{X}_1*(\mathcal{A}*\mathcal{X}_2*\mathcal{A})*\mathcal{X}_1=\mathcal{X}_1
  *(\mathcal{A}*\mathcal{X}_2)^H*(\mathcal{A}*\mathcal{X}_1)^H \\
   &=& \mathcal{X}_1*(\mathcal{A}*\mathcal{X}_1*\mathcal{A}*\mathcal{X}_2)^H=\mathcal{X}_1*(\mathcal{A}*\mathcal{X}_2)^H \\
   &=& \mathcal{X}_1*\mathcal{A}*\mathcal{X}_2\\
   &=&  \mathcal{X}_1*(\mathcal{A}*\mathcal{X}_2*\mathcal{A})*\mathcal{X}_2=(\mathcal{X}_1*\mathcal{A})^H*(\mathcal{X}_2*\mathcal{A})^H*\mathcal{X}_2\\
   &=&  (\mathcal{X}_2*\mathcal{A}*\mathcal{X}_1*\mathcal{A})^H*\mathcal{X}_2=(\mathcal{X}_2*\mathcal{A})^H*\mathcal{X}_2\\
   &=& \mathcal{X}_2*\mathcal{A}*\mathcal{X}_2=\mathcal{X}_2.
\end{eqnarray*}
Therefore, the Moore-Penrose inverse of $\mathcal{A}$ is unique. $\Box$


In fact, the tensor $\mathcal{R}$ obtained in the proof of Theorem \ref{t2} is the Moore-Penrose inverse of the tensor $\mathcal{S}$. So, the following is straightforward.

\begin{corollary}\label{c0}
Let $\mathcal{A}\in\mathbb{H}^{n_1\times n_2\times n_3}$ and $\mathcal{A}=\mathcal{U}*\mathcal{S}*\mathcal{V}^H$ be the SVD of $\mathcal{A}$.
Then,
\begin{equation*}
  \mathcal{A}^\dag=\mathcal{V}*\mathcal{S}^\dag*\mathcal{U}^H.
\end{equation*}

\end{corollary}

By using Theorem \ref{vt} and Corollary \ref{c0}, we can prove the following context.

\begin{theorem}\label{pp}
Let $\mathcal{A}\in\mathbb{H}^{n_1\times n_2\times n_3}$. Then,
\begin{itemize}
\item [] $(i)$  $(\mathcal{A}^{(1)})^H\in\mathcal{A}^H\{1\}$. In particular, $(\mathcal{A}^H)^\dag=(\mathcal{A}^\dag)^H$.
  \item [] $(ii)$ $(\mathcal{A}^H*\mathcal{A})^\dag=\mathcal{A}^\dag*(\mathcal{A}^H)^\dag$,
      $(\mathcal{A}*\mathcal{A}^H)^\dag=(\mathcal{A}^H)^\dag*\mathcal{A}^\dag$.
  \item [] $(iii)$ $\mathcal{A}^H=\mathcal{A}^H*\mathcal{A}*\mathcal{A}^\dag
=\mathcal{A}^\dag*\mathcal{A}*\mathcal{A}^H$.
  \item [] $(iv)$ $\mathcal{A}^\dag=(\mathcal{A}^H*\mathcal{A})^\dag*\mathcal{A}^H
=\mathcal{A}^H*(\mathcal{A}*\mathcal{A}^H)^\dag$.
\end{itemize}
\end{theorem}
{\sc Proof:} (i) It follows by the definition of the $\{1\}$-inverse and the Moore-Penrose inverse of a quaternion tensor.

(ii) Let $\mathcal{A}=\mathcal{U}*\mathcal{S}*\mathcal{V}^H$ be the SVD of $\mathcal{A}$. Then, by Corollary \ref{c0}, $\mathcal{A}^\dag=\mathcal{V}*\mathcal{S}^\dag*\mathcal{U}^H$. Meanwhile, we have $$\mathcal{A}^H=\mathcal{V}*\mathcal{S}^H*\mathcal{U}^H \  \text{and} \ (\mathcal{A}^H)^\dag=\mathcal{U}*(\mathcal{S}^H)^\dag*\mathcal{V}^H.$$
Then,
\begin{eqnarray*}
  (\mathcal{A}^H*\mathcal{A})^\dag &=& (\mathcal{V}*\mathcal{S}^H*\mathcal{U}^H
*\mathcal{U}*\mathcal{S}*\mathcal{V}^H)^\dag=
\mathcal{V}*(\mathcal{S}^H*\mathcal{S})^\dag*\mathcal{V}^H \\
   &=& \mathcal{V}*\mathcal{S}^\dag*(\mathcal{S}^H)^\dag*\mathcal{V}^H
   =  \mathcal{V}*\mathcal{S}^\dag*\mathcal{U}^H*\mathcal{U}*(\mathcal{S}^H)^\dag*\mathcal{V}^H\\
  &=& \mathcal{A}^\dag*(\mathcal{A}^H)^\dag.
\end{eqnarray*}
Similarly, we have $(\mathcal{A}*\mathcal{A}^H)^\dag=(\mathcal{A}^H)^\dag*\mathcal{A}^\dag$.

(iii) Notice that
\begin{eqnarray*}
  \mathcal{A}^H*\mathcal{A}*\mathcal{A}^\dag &=& \mathcal{V}*\mathcal{S}^H*\mathcal{U}^H*
   \mathcal{U}*\mathcal{S}*\mathcal{V}^H*\mathcal{V}*\mathcal{S}^\dag*\mathcal{U}^H\\
   &=& \mathcal{V}*\mathcal{S}^H*\mathcal{S}*\mathcal{S}^\dag*\mathcal{U}^H=
   \mathcal{V}*\mathcal{S}^H*(\mathcal{S}^\dag)^H*\mathcal{S}^H*\mathcal{U}^H\\
   &=& \mathcal{V}*\mathcal{S}^H*\mathcal{U}^H=\mathcal{A}^H.
\end{eqnarray*}
Similarly, we have $\mathcal{A}^H=\mathcal{A}^\dag*\mathcal{A}*\mathcal{A}^H$.

(iv) By using (ii), we have
\begin{equation*}
  (\mathcal{A}^H*\mathcal{A})^\dag*\mathcal{A}^H=\mathcal{A}^\dag*(\mathcal{A}^H)^\dag*\mathcal{A}^H
  =\mathcal{A}^\dag*(\mathcal{A}*\mathcal{A}^\dag)^H=\mathcal{A}^\dag*\mathcal{A}*\mathcal{A}^\dag
  =\mathcal{A}^\dag.
\end{equation*}
Analogously, $\mathcal{A}^\dag=\mathcal{A}^H*(\mathcal{A}*\mathcal{A}^H)^\dag$. $\Box$

By using the same technical of complex matrices, we can prove the following theorem. The reader can consult \cite{B-G}.
\begin{theorem}\label{tm4}
Let $\mathcal{A}\in\mathbb{H}^{n_1\times l\times n_3}$, $\mathcal{B}\in\mathbb{H}^{m\times k\times n_3}$, $\mathcal{C}\in\mathbb{H}^{n_1\times k\times n_3}$ and $\mathcal{X}\in\mathbb{H}^{l\times m\times n_3}$. Then, the tensor equation
\begin{equation*}\label{vv}
  \mathcal{A}*\mathcal{X}*\mathcal{B}=\mathcal{C}
\end{equation*}
is consistent if and only if exist $\mathcal{A}^{(1)}\in\mathcal{A}\{1\}$, $\mathcal{B}^{(1)}\in\mathcal{B}\{1\}$ such that
\begin{equation*}\label{vv2}
  \mathcal{A}*\mathcal{A}^{(1)}*\mathcal{C}*\mathcal{B}^{(1)}*\mathcal{B}=\mathcal{C}
\end{equation*}
in which case the general solution is
\begin{equation}\label{vv3}
  \mathcal{X}=\mathcal{A}^{(1)}*\mathcal{C}*\mathcal{B}^{(1)}+\mathcal{W}-
  \mathcal{A}^{(1)}*\mathcal{A}*\mathcal{W}*\mathcal{B}*\mathcal{B}^{(1)}
\end{equation}
for arbitrary $\mathcal{W}\in\mathbb{H}^{l\times m\times n_3}$.
\end{theorem}

\begin{theorem}\label{wa}
Let $\ac \in
\mathbb{H}^{n_1\times n_2\times n_3}$. Then, the set $\mathcal{A}\{1, 3\}$ consists of all solutions $\mathcal{X}$ of
\begin{equation*}
  \mathcal{A}*\mathcal{X}
=\mathcal{A}*\mathcal{A}^{(1,3)}.
 \end{equation*}
\end{theorem}

{\sc Proof:} Since
$$\mathcal{A}*\mathcal{X}*\mathcal{A}=\mathcal{A}*\mathcal{A}^{(1,3)}*\mathcal{A}=\mathcal{A},$$
and $(\mathcal{A}*\mathcal{X})^H=\mathcal{A}*\mathcal{X}$ due to $(\mathcal{A}*\mathcal{A}^{(1,3)})^H=\mathcal{A}*\mathcal{A}^{(1,3)}$, then $\mathcal{X}\in \mathcal{A}\{1,3\}$.

If $\mathcal{X}\in \mathcal{A}\{1,3\}$, then
$$\mathcal{A}*\mathcal{A}^{(1,3)}=\mathcal{A}*\mathcal{X}*\mathcal{A}*\mathcal{A}^{(1,3)}
=\mathcal{X}^H*\mathcal{A}^H*(\mathcal{A}^{(1,3)})^H*\mathcal{A}^H=
\mathcal{X}^H*\mathcal{A}^H=\mathcal{A}*\mathcal{X},$$
where we have used Theorem \ref{pp} (i). $\Box$

The following theorem is obtained in a manner analogous to Theorem \ref{wa}.

\begin{theorem}
Let $\ac \in
\mathbb{H}^{n_1\times n_2\times n_3}$. Then, the set $\mathcal{A}\{1, 4\}$ consists of all solutions $\mathcal{X}$ of
\begin{equation*}
\mathcal{X}*\mathcal{A}
=\mathcal{A}^{(1,4)}*\mathcal{A}.
 \end{equation*}
\end{theorem}

\begin{theorem}
Let $\ac \in
\mathbb{H}^{n_1\times n_2\times n_3}$. Then, the following statements are true. \begin{itemize}
  \item [] $(i)$ $\mathcal{A}\{1\}=\{\mathcal{A}^{(1)}+\mathcal{Z}-
\mathcal{A}^{(1)}*\mathcal{A}*\mathcal{Z}*\mathcal{A}*\mathcal{A}^{(1)}: \mathcal{Z} \in\mathbb{H}^{n_2\times n_1\times n_3}\}$.
  \item [] $(ii)$ $\mathcal{A}\{1,3\}=\{\mathcal{A}^{(1,3)}+
(\mathcal{I}-\mathcal{A}^{(1,3)}*\mathcal{A})*\mathcal{Z}: \mathcal{Z}\in\mathbb{H}^{n_2\times n_1\times n_3}\}$.
  \item [] $(iii)$ $\mathcal{A}\{1,4\}=\{\mathcal{A}^{(1,4)}+
\mathcal{Z}*(\mathcal{I}-\mathcal{A}*\mathcal{A}^{(1,4)}): \mathcal{Z}\in\mathbb{H}^{n_2\times n_1\times n_3}\}$.
\end{itemize}
\end{theorem}
{\sc Proof:} (i) By Theorem \ref{tm4}, the general solution of $\mathcal{A}*\mathcal{X}*\mathcal{A}=\mathcal{A}$ is
$$\mathcal{X}=\mathcal{A}^{(1)}*\mathcal{A}*\mathcal{A}^{(1)}+\mathcal{W}-
  \mathcal{A}^{(1)}*\mathcal{A}*\mathcal{W}*\mathcal{A}*\mathcal{A}^{(1)}, \ \text{where} \ \mathcal{W}\in\mathbb{H}^{n_2\times n_1\times n_3}.$$
Taking $\mathcal{W}=\mathcal{A}^{(1)}+\mathcal{Z}$ in the
set of solutions of $\mathcal{A}*\mathcal{X}*\mathcal{A}=\mathcal{A}$, we get
$$\mathcal{A}\{1\}=\{\mathcal{A}^{(1)}+\mathcal{Z}-
\mathcal{A}^{(1)}*\mathcal{A}*\mathcal{Z}*\mathcal{A}*\mathcal{A}^{(1)}: \mathcal{Z} \in\mathbb{H}^{n_2\times n_1\times n_3}\}.$$

(ii) By Theorem \ref{tm4}, the general solution of $\mathcal{A}*\mathcal{X}
=\mathcal{A}*\mathcal{A}^{(1,3)}$ is
$$\mathcal{X}=\mathcal{A}^{(1)}*\mathcal{A}*\mathcal{A}^{(1,3)}+\mathcal{W}-
  \mathcal{A}^{(1)}*\mathcal{A}*\mathcal{W}, \ \text{where} \ \mathcal{W}\in
\mathbb{H}^{n_2\times n_1\times n_3}.$$
Substituting $\mathcal{Z}+\mathcal{A}^{(1,3)}$ for $\mathcal{W}$ gives
$\mathcal{X}=\mathcal{A}^{(1,3)}+
(\mathcal{I}-\mathcal{A}^{(1,3)}*\mathcal{A})*\mathcal{Z}$, that is
$$\mathcal{A}\{1,3\}=\{\mathcal{A}^{(1,3)}+
(\mathcal{I}-\mathcal{A}^{(1,3)}*\mathcal{A})*\mathcal{Z}: \mathcal{Z}\in\mathbb{H}^{n_2\times n_1\times n_3}\}.$$

(iii) Similar as (ii).

\section{An Algorithm for Computing the Moore-Penrose Inverse of a Quaternion Tensor}

By refer to Algorithm 2 in \cite{DK}, we propose the following Algorithm \ref{algo:sylvester1} for computing the Moore-Penrose inverse of a quaternion tensor. Before that, we declare that $\fft(\cdot)$ and $\ifft(\cdot)$ are Matlab functions, which implement the fast Fourier transform and the inverse fast Fourier transform of a matrix, respectively.

\begin{algorithm}[H]
\caption{\textsc{Compute the Moore-Penrose inverse of a quaternion tensor $\mathcal{A}$}}\label{algo:sylvester1}
\KwIn{$n_1\times n_2\times n_3$ quaternion tensor $\mathcal{A}$}
\KwOut{$n_2\times n_1\times n_3$ quaternion tensor $\mathcal{X}$}
\begin{enumerate}
\addtolength{\itemsep}{-0.8\parsep minus 0.8\parsep}

\item $i=3$, $\mathcal{D} = \fft(\mathcal{A},[ ~ ],i)$;

\item  for $i=1,\ldots, n_3$

\qquad  $[U_i,S_i,V_i]$=svd$(\mathcal{D}(:,:,i))$;

\qquad  $R_i=$pinv$(S_i)$; where pinv$(S_i)$ is the Moore-Penrose inverse of
$S_i$

\qquad $\mathcal{G}(:,:,i)=V_iR_iU^H_i$;

end

\item $i=3$, $\mathcal{X}=\ifft(\mathcal{G},[ ~ ],i)$.
\end{enumerate}
\end{algorithm}

 The strategy of this algorithm is using $\fft(\cdot)$ to some objects and then calculate the Moore-Penrose inverse of each result matrix from $\fft(\mathcal{A})$. Finally, employing $\ifft(\cdot)$ to $\mathcal{G}(:,:,i)$ as in the Algorithm  to get the Moore-Penrose inverse of $\mathcal{A}$. Next, we will test the construct Algorithm \ref{algo:sylvester1} by using the following example.

\begin{myexample}
Let $\mathcal{A}$ be a $2\times3\times4$ quaternion tensor with the following form:
\begin{equation*}
  \mathcal{A}(:,:,1)=\begin{bmatrix}
                             {\bf i} &{\bf j} &2{\bf i}-{\bf j}\\
                             {\bf k} &2+{\bf i} &3
                          \end{bmatrix}, \qquad \mathcal{A}(:,:,2)=\begin{bmatrix}
                               1+{\bf i} &2 &1+{\bf i}+{\bf j}+{\bf k}\\
                               {\bf k} &{\bf j}+3{\bf k} &1+{\bf k}\\
                          \end{bmatrix},
\end{equation*}

\begin{equation*}
  \mathcal{A}(:,:,3)=\begin{bmatrix}
                              1+3{\bf j}  &  3{\bf j}-{\bf k}  & 5 \\
   2+{\bf i}+{\bf j}+2{\bf k}  &  {\bf i}+{\bf j} &  {\bf k}  \\
                          \end{bmatrix}, \qquad  \mathcal{A}(:,:,4)=\begin{bmatrix}
                             2{\bf k} &  {\bf i}-{\bf j}-{\bf k}  &  -2{\bf j} \\
   3{\bf j}-{\bf k}  & {\bf i}+{\bf k}  &  3{\bf i}-{\bf j} \\
                          \end{bmatrix}.
\end{equation*}
Implement Algorithm \ref{algo:sylvester1} on $\mathcal{A}$, we have

\begin{equation*}
  \mathcal{A}^\dag(:,:,1)=\begin{bmatrix}
   0.0588-0.0539{\bf i}-0.0196{\bf j}-0.0049{\bf k}  &   -0.0245-0.0196{\bf i}-0.0049{\bf j}+0.0294{\bf k}   \\
   -0.0221+0.0025{\bf i}+0.0368{\bf j}-0.0662{\bf k}  &   0.0270+0.0074{\bf i}v+0.0074{\bf j}-0.0613{\bf k}   \\
   -0.0735-0.0588{\bf i}-0.0686{\bf j}+0.0196{\bf k}  &   0.0490-0.0980{\bf i}-0.0392{\bf j}+0.0833{\bf k}   \\
                          \end{bmatrix},
\end{equation*}

\begin{equation*}
  \mathcal{A}^\dag(:,:,2)=\begin{bmatrix}
 0.0526-0.0395{\bf i}-0.0175{\bf j}-0.0614{\bf k}  &   -0.0175-0.0219{\bf i}-0.0263{\bf j}+0.0175{\bf k}   \\
 -0.0461-0.0154{\bf i}+0.0329{\bf j}-0.0241{\bf k}  &   -0.0811-0.0066{\bf i}-0.0110{\bf j}+0.0461{\bf k}   \\
 0.0197-0.0110{\bf i}+0.0241{\bf j}+0.0855{\bf k}  &   0.1162-0.0636{\bf i}-0.0022{\bf j}-0.0373{\bf k}   \\
                          \end{bmatrix},
\end{equation*}

\begin{equation*}
  \mathcal{A}^\dag(:,:,3)=\begin{bmatrix}
 0.0238-0.0087{\bf i}-0.0389{\bf j}-0.0325{\bf k}  &   0.0198+0.0246{\bf i}+0.0167{\bf j}+0.0381{\bf k}   \\
 0.0329-0.0115{\bf i}-0.0329{\bf j}-0.0202{\bf k}  &   -0.0448-0.0044{\bf i}+0.0179{\bf j}-0.0258{\bf k}   \\
 0.0111-0.0214{\bf i}+0.0349{\bf j}-0.0246{\bf k}  &   -0.0595+0.0333{\bf i}-0.0071{\bf j}+0.0651{\bf k}   \\
                          \end{bmatrix},
\end{equation*}

\begin{equation*}
  \mathcal{A}^\dag(:,:,4)=\begin{bmatrix}
 0.0452-0.0393{\bf i}+0.0357{\bf j}+0.0083{\bf k}  &   0.0310+0.0202{\bf i}+0.0548{\bf j}-0.0464{\bf k}   \\
 0.0286+0.0512{\bf i}-0.0679{\bf j}-0.0643{\bf k}  &   -0.0262+0.0607{\bf i}-0.0298{\bf j}-0.0238{\bf k}   \\
 0.0226+0.0417{\bf i}-0.0107{\bf j}+0.0131{\bf k}  &   -0.0083+0.0012{\bf i}+0.0679{\bf j}-0.0512{\bf k}   \\
                          \end{bmatrix}.
\end{equation*}

\end{myexample}

\section{The Drazin Inverse of Quaternion Tensors}\label{three}

Let us recall that the index of a matrix $A$ is the smallest nonnegative integer $k$ such that
$\rk(A^k) = \rk(A^{k+1})$. It is denoted by $\ind(A)$.

\begin{definition} 
Let $\ac \in \mathbb{H}^{n_1\times n_2\times n_3}$ and
\begin{eqnarray*} 
\DFT(\ccc(\unfold(\mathcal{A}))) =
\begin{bmatrix} A_1 & & \\ & \ddots & \\ & & A_{n_3} \end{bmatrix}.
\end{eqnarray*}
The ${n_3}$-tuple
\begin{equation*}
\big(\ind(A_1), \ind(A_2),\ldots, \ind(A_{n_3})\big)
\end{equation*}
is called the {\bf{multi-index}} of the quaternion tensor $\ac$.
\end{definition}

\begin{definition} \label{d52}
Let $\ac, \xc \in \mathbb{H}^{n_1\times n_1\times n_3}$. Suppose
that $\ac$ and $\xc$ are expressed as follows:
\begin{eqnarray*} 
\DFT(\ccc(\unfold(\ac))) =
\begin{bmatrix} A_1 & & \\ & \ddots & \\ & & A_{n_3} \end{bmatrix}
\quad \text{and} \quad
\DFT(\ccc(\unfold(\xc))) =
\begin{bmatrix} X_1 & & \\ & \ddots & \\ &  &  X_{n_3} \end{bmatrix}.
\end{eqnarray*}
If $X_i$ is the Drazin inverse of $A_i$ for $i=1, \ldots, {n_3}$, i.e.,
\begin{equation}\label{mpk3}
A_i X_i A_i^{l_i} = A_i^{l_i}, \quad X_iA_iX_i = X_i, \quad A_iX_i = X_i A_i, \qquad i=1,\cdots,{n_3},
\end{equation}
where $l_i = \ind(A_i)$, then $\xc$ is called the {\bf{Drazin inverse}} of the quaternion
tensor $\ac$ and is denoted by $\ac^D$.
\end{definition}

\begin{lemma}
Let $\ac \in
\mathbb{H}^{n_1\times n_1\times n_3}$ and
\begin{eqnarray*} 
\DFT(\ccc(\unfold(\mathcal{A}))) =
\begin{bmatrix} A_1 & & \\ & \ddots & \\ & & A_{n_3}\end{bmatrix}.
\end{eqnarray*}
If the multi-index of $\mathcal{A}$ is $(l_1,l_2,...,l_{n_3})$, then the Drazin inverse of $\mathcal{A}$  exists and is unique.
\end{lemma}
{\sc Proof:} Since $\ind(A_i) = l_i$ for $1 \leq i \leq {n_3}$, the matrices $A_1, \ldots, A_{n_3}$ are Drazin
invertible. Let $X_i = A_i^D$ for $1 \leq i \leq {n_3}$ and define the quaternion tensor $\mathcal{X}$
as in Definition \ref{d52}. It is trivial to see that $\xc$ is the Drazin inverse of $\ac$.

Suppose both quaternion tensors $\mathcal{X}$, and $\mathcal{Y}$ satisfy the Definition \ref{d52}. Let
\begin{eqnarray*} 
\textbf{DFT}(\ccc(\textbf{Unfold}(\xc))) =
\begin{bmatrix} X_1 & & \\ & \ddots & \\ & & X_{n_3} \end{bmatrix}
\quad \text{and} \quad
\textbf{DFT}(\ccc(\textbf{Unfold}(\mathcal{Y}))) =
\begin{bmatrix} Y_1 & & \\ & \ddots & \\ &  & Y_{n_3} \end{bmatrix}.
\end{eqnarray*}
It follows $X_i = A_i^D$, $i=1,2,...,n_3$ and $Y_i = A_i^D$, $i=1,2,...,n_3$ satisfy (\ref{mpk3}). Hence, $\mathcal{X}$ and $\mathcal{Y}$ coincides since $X_i$ and $Y_i$ are the same. $\Box$

\begin{definition} 
Let $\ac \in \mathbb{H}^{n_1\times n_1\times n_3}$. If there exists a quaternion tensor
$\xc\in\mathbb{H}^{n_1\times n_1\times n_3}$ such that
\begin{equation*} 
\ac*\xc*\ac = \ac, \qquad \xc*\ac*\xc = \xc, \qquad \ac*\xc = \xc*\ac
\end{equation*}
then $\xc$ is called the {\bf{group inverse}} of the quaternion tensor $\ac$ and is denoted by $\ac^\#$.
\end{definition}

As for matrices, it can be proved easily that the group inverse of quaternion tenosrs, if exists, is unique.

\begin{lemma} 
Let $\ac \in \mathbb{H}^{n_1\times n_1\times n_3}$ and
\begin{eqnarray*} 
\DFT(\ccc(\unfold(\mathcal{A}))) =
\begin{bmatrix} A_1 & & \\ & \ddots & \\ & & A_{n_3}\end{bmatrix}.
\end{eqnarray*}
If
$\ind(A_i) = 1$ for $i=1, \ldots,{n_3}$, then the group inverse of the quaternion tensor $\ac$ exists and unique.
\end{lemma}

In the following, we will give some characterizations of the Drazin inverse  of the quaternion
tensor.

\begin{theorem}\label{JCF}{\rm\cite{Huang}} (Quaternion matrix Jordan Canonical Form) Let $A\in{H^{n\times n}}$. Then there exists
an invertible quaternion matrices $P\in{H^{n\times n}}$
such that
\begin{equation}\label{PAP}
  P^{-1}AP=\begin{bmatrix}
  J_{n_1}(\lambda_1) & 0 &\ldots&0  \\
  0 & J_{n_2}(\lambda_2) &\ldots&0   \\
  \vdots & \vdots &      & \vdots \\
   0 & 0 & \ldots  & 0 \\
  0 & 0 &\ldots&J_{n_s}(\lambda_s)   \\
\end{bmatrix}=J,
\end{equation}
where $\lambda_k=a_k+b^2_k{\bf{i}}\in\mathbb{C}$, $a_k, b^2_k\in\mathbb{R}$, $k=1,2,...,s$, $\sum\limits_{k=1}^{s}n_k=n$. The $\lambda_k$, $k=1,2,...,s$ are all right eigenvalues of $A$ which are not necessarily distinct. The $J$ is uniquely determined by $A$ up to the order of matrices $J_{n_k}(\lambda_k)$ in (\ref{PAP}), where
$$J_{n_k}(\lambda_k)=\left[
                       \begin{array}{cccc}
                         \lambda_k & 1 &  &  \\
                          & \lambda_k & \ddots &  \\
                          &  & \ddots & 1 \\
                          &  &  & \lambda_k \\
                       \end{array}
                     \right]\in{C}^{n_k\times n_k}.
$$

\end{theorem}

\begin{theorem}
Let $\ac \in \mathbb{H}^{n_1\times n_1\times n_3}$,
then there
exists an invertible quaternion tensor $\mathcal{P} \in \mathbb{H}^{n_1\times n_1\times n_3}$ such that
\begin{eqnarray*} 
\mathcal{A}=\mathcal{P}^{-1}*\mathcal{J}*\mathcal{P},
\end{eqnarray*}
where $\mathcal{J}=\fold(\uncirc(\IDFT(\begin{bmatrix}
                     J_1 &   &   \\
                       & \ddots   &   \\
                       &   &  J_{n_3} \\
                   \end{bmatrix})))$, $J_i$, $i=1,2,...,n_3$  have the forms as $J$ defined in (\ref{PAP}).
\end{theorem}
{\sc Proof:} Let \begin{eqnarray*} 
\textbf{DFT}(\text{circ}(\unfold(\mathcal{A}))) =
\begin{bmatrix} A_1 & & \\ & \ddots & \\ & & A_{n_3} \end{bmatrix}.
\end{eqnarray*}
By Theorem \ref{JCF}, each $A_i$, $i=1,2,...,{n_3}$ have the decomposition $A_i=P_i^{-1}J_iP_i$, where $P_i$ are invertible quaternion matrices, $J_i$ are defined as in (\ref{PAP}). Thus, we have
\begin{equation}\label{PAP2}
  \begin{bmatrix}
                     A_1 &   &   \\
                       & \ddots   &   \\
                       &   &  A_{n_3} \\
                   \end{bmatrix}=\begin{bmatrix}
                     P^{-1}_1 &   &   \\
                       & \ddots   &   \\
                       &   &  P^{-1}_{n_3} \\
                   \end{bmatrix}\begin{bmatrix}
                     J_1 &   &   \\
                       & \ddots   &   \\
                       &   &  J_{n_3} \\
                   \end{bmatrix}\begin{bmatrix}
                     P_1 &   &   \\
                       & \ddots   &   \\
                       &   &  P_{n_3} \\
                   \end{bmatrix}.
\end{equation}
Implementing $\fold(\text{uncirc}(\textbf{IDFT}))(\cdot)$ on both sides of the equation above, we get $\mathcal{A}=\mathcal{P}^{-1}*\mathcal{J}*\mathcal{P}$. $\Box$

\begin{theorem}
Let $\ac \in \mathbb{H}^{n_1\times n_1\times n_3}$ and the multi-index of $\mathcal{A}$ is $(l_1,l_2,...,l_{n_3})$,
then
\begin{eqnarray*} 
\mathcal{A}^D=\mathcal{P}^{-1}*\mathcal{J}^D*\mathcal{P},
\end{eqnarray*}
where $\mathcal{P} \in \mathbb{H}^{n_1\times n_1\times n_3}$ is  an invertible quaternion tensor and $$\mathcal{J}^D=\fold(\uncirc(\IDFT(\begin{bmatrix}
                     J^D_1 &   &   \\
                       & \ddots   &   \\
                       &   &  J^D_{n_3} \\
                   \end{bmatrix}))),$$
                   $J_i$, $i=1,2,...,{n_3}$  have the forms as $J$ defined in (\ref{PAP}).
\end{theorem}
{\sc Proof:}
Let \begin{eqnarray*} 
\DFT(\ccc(\unfold(\mathcal{A}))) =
\begin{bmatrix} A_1 & & \\
 & \ddots & \\
 & & A_{n_3} \end{bmatrix}.
\end{eqnarray*}
Then, by using Theorem \ref{JCF}, one has
\begin{eqnarray*}
\begin{bmatrix} A_1 & & \\
 & \ddots & \\
 & & A_{n_3} \end{bmatrix}&=&\begin{bmatrix} P^{-1}_1J_1P_1 & & \\
 & \ddots & \\
 & & P^{-1}_{n_3} J_{n_3} P_{n_3} \end{bmatrix}\\
 &=&\begin{bmatrix}
                     P^{-1}_1 &   &   \\
                       & \ddots   &   \\
                       &   &  P^{-1}_{n_3} \\
                   \end{bmatrix}\begin{bmatrix}
                     J_1 &   &   \\
                       & \ddots   &   \\
                       &   &  J_{n_3} \\
                   \end{bmatrix}\begin{bmatrix}
                     P_1 &   &   \\
                       & \ddots   &   \\
                       &   &  P_{n_3} \\
                   \end{bmatrix}.
\end{eqnarray*}
 It is easy to check
\begin{eqnarray*}
\begin{bmatrix} A^D_1 & & \\
 & \ddots & \\
 & & A^D_{n_3} \end{bmatrix}=\begin{bmatrix}
                     P^{-1}_1 &   &   \\
                       & \ddots   &   \\
                       &   &  P^{-1}_{n_3} \\
                   \end{bmatrix}\begin{bmatrix}
                     J^D_1 &   &   \\
                       & \ddots   &   \\
                       &   &  J^D_{n_3} \\
                   \end{bmatrix}\begin{bmatrix}
                     P_1 &   &   \\
                       & \ddots   &   \\
                       &   &  P_{n_3} \\
                   \end{bmatrix}.
\end{eqnarray*}
Hence, implementing $\fold(\text{uncirc}(\textbf{IDFT}))(\cdot)$ on both sides of the equation above, we get $\mathcal{A}^D=\mathcal{P}^{-1}*\mathcal{J}^D*\mathcal{P}$. $\Box$

In the following, we will give another representation of the Drazin inverse  by using the core-nilpotent of the quaternion tensors.

\begin{definition}
Let $\ac \in \mathbb{H}^{n_1\times n_1\times n_3}$.
The product
$$\mathcal{C}_{\mathcal{A}} = \mathcal{A}^2
*\mathcal{A}^D =\fold(\uncirc(\IDFT(\begin{bmatrix} A^2_1A^D_1 & & \\ & \ddots & \\ &  & A^2_{n_3} A^D_{n_3} \end{bmatrix})))
$$
is called the {\bf{core}} of the quaternion tensor $\mathcal{A}$.
\end{definition}


\begin{definition}
Let $\ac \in \mathbb{H}^{n_1\times n_1\times n_3}$ and $\mathcal{C}_{\mathcal{A}}\in \mathbb{H}^{n_1\times n_1\times n_3}$ is the core of the quaternion tensor $\mathcal{A}$. Suppose $\ac$  and $\mathcal{C}_{\mathcal{A}}$ have the following forms
\begin{eqnarray*} 
\fold(\uncirc(\IDFT(\ac))) =
\begin{bmatrix} A_1 & & \\ & \ddots & \\ & & A_\rho \end{bmatrix},  \quad
\fold(\uncirc(\IDFT(\mathcal{C}_{\mathcal{A}}))) =
\begin{bmatrix} A^2_1A^D_1 & & \\ & \ddots & \\ &  & A^2_\rho A^D_\rho \end{bmatrix}.
\end{eqnarray*}
Then,
$$\mathcal{N}_{\mathcal{A}} = \mathcal{A}-\mathcal{C}_{\mathcal{A}}=\fold(\uncirc(\IDFT(\begin{bmatrix}
                     A_1-A^2_1A^D_1 &   &   \\
                       & \ddots   &   \\
                       &   &  A_\rho-A^2_\rho A^D_\rho \\
                   \end{bmatrix})))$$
is called a {\bf{nilpotent}} quaternion tensor.
\end{definition}

The following inspired by \cite{Wang} is easy to check by using the technique of complex matrices.
\begin{theorem}\label{taa} Let $A\in{H^{n\times n}}$, $\ind(A)=k$, and $A=C_A+N_A$  is the core-nilpotent decomposition of $A$. Then, there exists an invertible quaternion matrices $P\in{H^{n\times n}}$ such that
\begin{eqnarray*}
A=P\begin{bmatrix}
  C & 0   \\
  0 & N   \\
\end{bmatrix}P^{-1},
\end{eqnarray*} where $C_A=P\begin{bmatrix}
  C & 0   \\
  0 & 0   \\
\end{bmatrix}P^{-1}$, $N_A=P\begin{bmatrix}
  0 & 0   \\
  0 & N   \\
\end{bmatrix}P^{-1}$, $C\in{H^{r\times r}}$, $N\in{H^{(n-r)\times (n-r)}}$ and then
\begin{eqnarray*}
A^D=P\begin{bmatrix}
  C^{-1} & 0   \\
  0 & 0   \\
\end{bmatrix}P^{-1}.
\end{eqnarray*}.

\end{theorem}

\begin{theorem}\label{t51}
Let $\ac \in \mathbb{H}^{n_1\times n_1\times n_3}$ and the multi-index of $\mathcal{A}$ is $(l_1,l_2,...,l_{n_3})$,
then
\begin{equation} \label{a}
\ac = \mathcal{P}*\Phi*\mathcal{P}^{-1},
\end{equation}
where $\mathcal{P} \in \mathbb{H}^{n_1\times n_1\times n_3}$ is  an invertible quaternion tensor,
\begin{equation*}
  \Phi=\fold(\uncirc(\IDFT(
\begin{bmatrix}
\mat{C_1}{0}{0}{N_1} &  & \\ & \ddots & \\ &  & \mat{C_{n_3}}{0}{0}{N_{n_3}}
\end{bmatrix}))),
\end{equation*}
$\mat{C_i}{0}{0}{N_i}= P^{-1}_i(C_{A_i}+N_{A_i})P_i$, $C_{A_i}$ and $N_{A_i}$ are the core and nilpotent part of $A_i$, $i=1,2,...,{n_3}$, respectively. Besides,
\begin{equation} \label{a}
\ac^D = \mathcal{P}*\Phi^D*\mathcal{P}^{-1},
\end{equation}
where
\begin{equation*}
  \Phi^D=\fold(\uncirc(\IDFT(
\begin{bmatrix}
\mat{C^{-1}_1}{0}{0}{0} &  & \\ & \ddots & \\ &  & \mat{C^{-1}_{n_3}}{0}{0}{0}
\end{bmatrix}))).
\end{equation*}
\end{theorem}
{\sc Proof:}
Let \begin{eqnarray*} 
\textbf{DFT}(\text{circ}(\unfold(\mathcal{A}))) =
\begin{bmatrix} A_1 & & \\
 & \ddots & \\
 & & A_{n_3} \end{bmatrix}.
\end{eqnarray*}
Then, by Theorem \ref{taa}, we have
\begin{eqnarray*}
\begin{bmatrix} A_1 & & \\
 & \ddots & \\
 & & A_{n_3} \end{bmatrix}&=&\begin{bmatrix}P_1\begin{bmatrix}
  C_1 & 0   \\
  0 & N_1   \\
\end{bmatrix}P_1^{-1} & & \\
 & \ddots & \\
 & & P_{n_3}\begin{bmatrix}
  C_{n_3} & 0   \\
  0 & N_{n_3}   \\
\end{bmatrix}P_{n_3}^{-1} \end{bmatrix}\\
&=&\begin{bmatrix}
                     P_1 &   &   \\
                       & \ddots   &   \\
                       &   &  P_{n_3} \\
                   \end{bmatrix}\begin{bmatrix}
                     \begin{bmatrix}
  C_1 & 0   \\
  0 & N_1   \\
\end{bmatrix} &   &   \\
                       & \ddots   &   \\
                       &   &  \begin{bmatrix}
  C_{n_3} & 0   \\
  0 & N_{n_3}   \\
\end{bmatrix} \\
                   \end{bmatrix}\begin{bmatrix}
                     P^{-1}_1 &   &   \\
                       & \ddots   &   \\
                       &   &  P^{-1}_{n_3} \\
                   \end{bmatrix}.
\end{eqnarray*}
Implementing $\fold(\text{uncirc}(\textbf{IDFT}))(\cdot)$ on the quaternion tensors of the both sides of equation, one has
$$\mathcal{A}=\mathcal{P}*\Phi*\mathcal{P}^{-1},$$
where \begin{equation*}
  \Phi=\fold(\uncirc(\IDFT(
\begin{bmatrix}
\mat{C_1}{0}{0}{N_1} &  & \\ & \ddots & \\ &  & \mat{C_{n_3}}{0}{0}{N_{n_3}}
\end{bmatrix}))).
\end{equation*}
Notice that by Theorem \ref{taa}, we have
\begin{eqnarray*}
\begin{bmatrix} A^D_1 & & \\
 & \ddots & \\
 & & A^D_{n_3} \end{bmatrix}&=&\begin{bmatrix}(P_1\begin{bmatrix}
  C_1 & 0   \\
  0 & N_1   \\
\end{bmatrix}P_1^{-1})^D & & \\
 & \ddots & \\
 & & (P_{n_3}\begin{bmatrix}
  C_{n_3} & 0   \\
  0 & N_{n_3}   \\
\end{bmatrix}P_{n_3}^{-1})^D \end{bmatrix}\\
&=&\begin{bmatrix}
                     P_1 &   &   \\
                       & \ddots   &   \\
                       &   &  P_{n_3} \\
                   \end{bmatrix}\begin{bmatrix}
                     (\begin{bmatrix}
  C_1 & 0   \\
  0 & N_1   \\
\end{bmatrix})^D &   &   \\
                       & \ddots   &   \\
                       &   &  (\begin{bmatrix}
  C_{n_3} & 0   \\
  0 & N_{n_3}   \\
\end{bmatrix})^D \\
                   \end{bmatrix}\begin{bmatrix}
                     P^{-1}_1 &   &   \\
                       & \ddots   &   \\
                       &   &  P^{-1}_{n_3} \\
                   \end{bmatrix}\\
&=&\begin{bmatrix}
                     P_1 &   &   \\
                       & \ddots   &   \\
                       &   &  P_{n_3} \\
                   \end{bmatrix}\begin{bmatrix}
                     \begin{bmatrix}
  C^{-1}_1 & 0   \\
  0 & 0   \\
\end{bmatrix} &   &   \\
                       & \ddots   &   \\
                       &   &  \begin{bmatrix}
  C^{-1}_{n_3} & 0   \\
  0 & 0   \\
\end{bmatrix} \\
                   \end{bmatrix}\begin{bmatrix}
                     P^{-1}_1 &   &   \\
                       & \ddots   &   \\
                       &   &  P^{-1}_{n_3} \\
                   \end{bmatrix}.
\end{eqnarray*}
Implementing $\fold(\text{uncirc}(\textbf{IDFT}))(\cdot)$ on both sides of the quaternion tensors, one has
\begin{equation} \label{a}
\ac^D=\mathcal{P}*\Phi^D*\mathcal{P}^{-1},
\end{equation}
where
\begin{equation*}
  \Phi^D=\fold(\uncirc(\IDFT(
\begin{bmatrix}
\mat{C^{-1}_1}{0}{0}{0} &  & \\ & \ddots & \\ &  & \mat{C^{-1}_{n_3}}{0}{0}{0}
\end{bmatrix}))).
\end{equation*}
$\Box$

Now, let us consider the partitioned form of a quaternion tensor. For an arbitrary quaternion tensor $\ac \in \mathbb{H}^{n_1\times n_2\times n_3}$, it can be written as
\begin{equation}\label{df}
\mathcal{A}=\left[
\begin{array}{cc}
\mathcal{A}_1 & \mathcal{A}_2 \\
\mathcal{A}_3 & \mathcal{A}_4 \\
\end{array}
\right],
\end{equation}
\text{where}  $\mathcal{A}_1\in\CC^{p\times q\times n_3}, \mathcal{A}_2\in\CC^{p\times (n_2-q)\times n_3}, \mathcal{A}_3\in\CC^{(n_1-p)\times q\times n_3},
\mathcal{A}_4\in\CC^{(n_1-p)\times(n_2-q)\times n_3}$.

Thus, for Theorem \ref{t51}, we can deduce if $rank(\mathbf{C}_i)=r$, $i=1,2,...,n_3$, then $\Phi$ can be written as
\begin{equation*}
  \Phi=\mat{\mathcal{C}}{\mathcal{O}}{\mathcal{O}}{\mathcal{N}}, \ \text{where} \
  \mathcal{C}\in\mathbb{H}^{r\times r\times n_3}, \ \mathcal{N}\in\mathbb{H}^{(n_1-r)\times(n_2-r)\times n_3}.
\end{equation*}
In this case,
\begin{equation*}
  \mathcal{A}=\mathcal{P}*\mat{\mathcal{C}}{\mathcal{O}}{\mathcal{O}}{\mathcal{N}}*\mathcal{P}^{-1}.
\end{equation*}
In addition,
\begin{equation*}
\ac^D=\mathcal{P}*\mat{\mathcal{C}^{-1}}{\mathcal{O}}{\mathcal{O}}{\mathcal{O}}*\mathcal{P}^{-1}.
\end{equation*}

\section{The Inverse Along Two Quaternion Tensors}
In this section, we will do some researches on the right/left inverse along two quaternion tensors.

\begin{definition}\label{dd53}
Let $\ac \in \mathbb{H}^{n_1\times n_2\times n_3}$ and
$\bc \in \mathbb{H}^{n_2\times l\times n_3}$,
$\cc \in \mathbb{H}^{k\times n_1\times n_3}$. If there exist tensors
$\zc \in \mathbb{H}^{n_2\times n_1\times n_3}$,
$\xc_1 \in \mathbb{H}^{l\times n_1\times n_3}$ and
$\yc_1 \in \mathbb{H}^{n_2\times k\times n_3}$ such that
\begin{equation} \label{matrizmary}
\zc*\ac*\bc = \bc, \qquad \cc*\ac*\zc = \cc, \qquad \zc = \bc*\xc_1=\yc_1*\cc,
\end{equation}
then $\zc$ is called the \textbf{right inverse along $\bc$, $\cc$ } and is denoted by $\mary{\ac_r}{\bc,\cc}$. In particular,
when $\ac \in{H}^{n_1\times n_2}$, $\zc\in{H}^{n_2\times n_1}$ is called the right inverse along the quaternion matrix $\bc,\cc$.
\end{definition}

\begin{definition} 
Let $\ac \in \mathbb{H}^{n_1\times n_2\times n_3}$ and
$\dc \in \mathbb{H}^{l\times n_1\times n_3}$,
$\ec \in \mathbb{H}^{n_2\times k\times n_3}$. If there exist tensors
$\zc \in \mathbb{H}^{n_2\times n_1\times n_3}$,
$\xc_2 \in \mathbb{H}^{n_2\times l\times n_3}$ and
$\yc_2 \in \mathbb{H}^{k\times n_1\times n_3}$ such that
\begin{equation} \label{matrizmary}
\dc*\ac*\zc = \dc, \qquad\zc*\ac*\ec = \ec, \qquad \zc = \xc_2*\dc=\ec*\yc_2,
\end{equation}
then $\zc$ is called the \textbf{left inverse along $\dc$, $\ec$ } and is denoted by $\mary{\ac_l}{\dc,\ec}$. In particular,
when $\ac \in{H}^{n_1\times n_2}$, $\zc\in{H}^{n_2\times n_1}$ is called the left inverse along the quaternion matrix $\dc,\ec$.
\end{definition}

In fact, the right/left inverse along two quaternion tensors if exists, is unique.
\begin{theorem}
Let $\ac \in \mathbb{H}^{n_1\times n_2\times n_3}$. Then the following statements are true.
\begin{itemize}
\item [] $(i)$ If $\ac$ is right invertible along $\mathcal{B},\mathcal{C}$, then it is unique.
\item [] $(ii)$ If $\ac$ is left invertible along $\mathcal{D},\mathcal{E}$, then it is unique.
\end{itemize}
\end{theorem}
{\sc Proof:} (i) Let $\zc_1, \zc_2\in\mathbb{H}^{n_2\times n_1\times n_3}$ be two right inverses of $\ac$ along
$\bc,\mathcal{C}$. Exist tensors $\xc_1,\xc_2, \yc_1, \yc_2$ of adequate size such that
\begin{equation*}
\zc_i*\ac*\bc = \bc, \qquad \cc*\ac*\zc_i = \cc, \qquad \zc_i = \bc*\xc_i=\yc_i*\cc,
\end{equation*}
for $i=1,2$. Now we have
\begin{equation*}
\zc_1 = \bc*\xc_1 = \zc_2*\ac*\bc*\xc_1 = \zc_2*\ac*\zc_1 = \yc_2*\cc*\ac*\zc_1 = \yc_2*\cc = \zc_2.
\end{equation*}
The proof is finished.

(ii) Similar to (i). $\Box$

The following definition is borrowed from \cite{SWC}. The readers can also refer to \cite{K}.

\begin{definition}{\rm\cite{SWC}}
For an arbitrary matrix $A\in H^{m\times n}$, we denote by
\begin{itemize}
\item [] $(i)$ $R_r(A)=\{y\in H^m : y=Ax, x\in H^n\}, \text{the column right space of A}$,

\item [] $(ii)$   $N_r(A)=\{x\in H^n : Ax=0\}, \text{the right null space of A}$,

\item [] $(iii)$  $R_l(A)=\{y\in H^n : y=xA, x\in H^m\}, \text{the row left space of A}$,

\item [] $(iv)$   $N_l(A)=\{x\in H^m : xA=0\}, \text{the left null space of A}$.
\end{itemize}
\end{definition}

In the paper \cite{SWC},  Song et al. defined ${A_r}^{(2)}_{T_1,S_1}$ and ${A_l}^{(2)}_{T_2,S_2}$ as follows.
\begin{definition}\label{ATS}{\rm\cite{SWC}} (i) An out inverse of a matrix $A\in H^{n_1\times n_2}$ with prescribed right range space $T_1$ and right null space $S_1$ is a solution of
the restricted matrix equation
\begin{equation}\label{ATS2}
  XAX=X,  \qquad R_r(X)=T_1,  \qquad N_r(X)=S_1
\end{equation}
and is denoted by $X={A_r}^{(2)}_{T_1,S_1}$.

(ii) An out inverse of a matrix $A\in H^{n_1\times n_2}$ with prescribed left range space $T_2$ and left null space $S_2$ is a solution of
the restricted matrix equation
\begin{equation}\label{ATS3}
  XAX=X,  \qquad R_l(X)=T_2,  \qquad N_l(X)=S_2
\end{equation}
and is denoted by $X={A_l}^{(2)}_{T_2,S_2}$.

\end{definition}

The right/left inverse along two quaternion matrices coincides with ${A_r}^{(2)}_{T_1,S_1}$/ ${A_l}^{(2)}_{T_2,S_2}$.
\begin{theorem}\label{TTT}
Let $A\in{H^{n_1\times n_2}}$, $B \in {H}^{n_2\times l}$,
$C \in {H}^{k\times n_1}$. The following statements are equivalent.
\begin{itemize}
\item [] $(i)$ $A$ is right invertible along $B,C$.
\item [] $(ii)$ The outer inverse ${A_r}^{(2)}_{T_1,S_1}$ exists.
\end{itemize}
\end{theorem}
{\sc Proof:} Suppose that statement (i) holds and let $Z=\mary{A_r}{B,C}$. Then,
\begin{equation*}
ZAZ=ZABX_1=BX_1=Z.
\end{equation*}
Notice that $Z=BX_1$ if and only if $R_r(Z) \subseteq R_r(B)$. $Z=Y_1C$ if and only if $N_r(Z) \supseteq N_r(C)$. In addition, $ZAB=B$, $CAZ=C$ is equivalent to $R_r(Z)\supseteq R_r(B)$,  $N_r(Z) \subseteq N_r(C)$. Hence, $R_r(Z)= R_r(B)$ and $N_r(Z) = N_r(C)$, that is $R_r(Z) = T_1$ and $N_r(Z) = S_1$, where $T_1$ and $S_1$ are defined in Definition \ref{ATS} (i). Hence, ${A_r}^{(2)}_{T_1,S_1}$ exists and ${A_r}^{(2)}_{T_1,S_1}=Z$.

Now suppose that statement (ii) holds and let $X={A_r}^{(2)}_{T_1,S_1}$.   In particular, $R_r(X)\subseteq R_r(B)$ and $N_r(C)\subseteq N_r(X)$. Then, there exists $X_1\in H^{l\times n_1}$ and $Y_1\in H^{n_2\times k}$
such that $X=BX_1$ and $X=Y_1C$. Also, $B$ and $C$ satisfy $B=XU$ and $C=VX$ for some $U\in H^{n_1\times l}$, $V\in H^{k\times n_2}$. This further implies
\begin{eqnarray*}
  B = XU= XAXU=XAB, \qquad C=VX=VXAX=CAX.
\end{eqnarray*}
Therefore, $X = \mary{A_r}{D_1,E_1}$. $\Box$

\begin{theorem}
Let $A\in{H^{n_1\times n_2}}$, $D \in {H}^{l\times n_1}$,
$E \in {H}^{n_2\times k}$. The following statements are equivalent.
\begin{itemize}
\item [] $(i)$ $A$ is left invertible along $D,E$.
\item [] $(ii)$ The outer inverse ${A_l}^{(2)}_{T_2,S_2}$ exists.
\end{itemize}
\end{theorem}
{\sc Proof:} Similar as Theorem \ref{TTT}. $\Box$


In the following, we will state some expressions for the right/left inverse along two quaternion matrices/tensors.

\begin{theorem}\label{t6} Let $A\in{H^{n_1\times n_2}}$, $B \in {H}^{n_2\times l}$,
$C \in {H}^{k\times n_1}$. If $A$ is right invertible along $B,C$, then
\begin{equation*}
  \mary{A_r}{B,C}=B(CAB)^\dag C.
\end{equation*}
\end{theorem}
{\sc Proof:} If $A$ is right invertible along $B,C$, then there exists $Z \in{H}^{n_2\times n_1}$ such that
\begin{equation} \label{m1}
ZAB= B, \qquad CAZ = C, \qquad Z = BX_1=Y_1C,
\end{equation}
where $X_1\in H^{l\times n_1}$ and $Y_1\in H^{n_2\times k}$.  Therefore, $C = CAZ = CABX_1$. In particular, $R_r(C) = R_r(CAB)$. Similarly, by $B=ZAB=Y_1CAB$, one has $N_r(B)=N_r(CAB)$. Then,
$$CAB(CAB)^\dag C=C, \qquad B(CAB)^\dag CAB=B.$$
Set $Z=B(CAB)^\dag C$, $X_1=(CAB)^\dag C$, $Y_1=B(CAB)^\dag$, it is easy to check $B(CAB)^\dag C$ is the solution of (\ref{m1}). Therefore, $\mary{A_r}{B,C}=B(CAB)^\dag C$. $\Box$




\begin{theorem}\label{tt4} Let $A\in{H^{n_1\times n_2}}$, $D \in {H}^{l\times n_1}$,
$E \in {H}^{n_2\times k}$. If $A$ is left invertible along $D,E$, then
\begin{equation*}
  \mary{A_l}{D,E}=E(DAE)^\dag D.
\end{equation*}
\end{theorem}
{\sc Proof:} Similar as Theorem \ref{t6}. $\Box$

\begin{theorem}\label{tt5}
Let $\ac \in \mathbb{H}^{n_1\times n_2\times n_3}$,
$\bc \in \mathbb{H}^{n_2\times l\times n_3}$,
$\cc \in \mathbb{H}^{k\times n_1\times n_3}$.
If $\mathcal{A}$ is right invertible along $\mathcal{B},\mathcal{C}$, then
\begin{equation*}
  \mary{\mathcal{A}_r}{\mathcal{B},\mathcal{C}}=\mathcal{B}*(\mathcal{C}*\mathcal{A}*\mathcal{B})^\dag *\mathcal{C},
\end{equation*}
where $(\mathcal{C}*\mathcal{A}*\mathcal{B})^\dag=\fold(\uncirc(\IDFT(\begin{bmatrix}
                     (C_1A_1B_1)^\dag &   &   \\
                       & \ddots   &   \\
                       &   &  (C_{n_3}A_{n_3}B_{n_3})^\dag \\
                   \end{bmatrix})))$.
\end{theorem}
{\sc Proof:} Let \begin{eqnarray*} 
\DFT(\ccc(\unfold(\mathcal{A}))) =
\begin{bmatrix} A_1 & & \\
 & \ddots & \\
 & & A_{n_3}\end{bmatrix}.
\end{eqnarray*}
Then, by Theorem \ref{t6}, we have
\begin{eqnarray*}
\begin{bmatrix} \mary{{A_1}_r}{B_1,C_1} & & \\
 & \ddots & \\
 & & \mary{{A_{n_3}}_r}{B_{n_3},C_{n_3}} \end{bmatrix}&=&\begin{bmatrix} B_1(C_1A_1B_1)^\dag C_1& & \\
 & \ddots & \\
 & & B_{n_3}(C_{n_3}A_{n_3}B_{n_3})^\dag C_{n_3}\end{bmatrix}\\
 &=&\begin{bmatrix}
                      B_1 &   &   \\
                       & \ddots   &   \\
                       &   &  B_{n_3} \\
                   \end{bmatrix}\begin{bmatrix}
                     (C_1A_1B_1)^\dag &   &   \\
                       & \ddots   &   \\
                       &   &  (C_{n_3}A_{n_3}B_{n_3})^\dag \\
                   \end{bmatrix}\begin{bmatrix}
                     C_1 &   &   \\
                       & \ddots   &   \\
                       &   &  C_{n_3} \\
                   \end{bmatrix}.
\end{eqnarray*}
Implementing $\fold(\text{uncirc}(\textbf{IDFT}))(\cdot)$ on the quaternion tensors of the both sides of the above equation, one has
$$\mary{\mathcal{A}_r}{\mathcal{B},\mathcal{C}}=\mathcal{B}*(\mathcal{C}*\mathcal{A}*\mathcal{B})^\dag *\mathcal{C},$$
where $(\mathcal{C}*\mathcal{A}*\mathcal{B})^\dag=\fold(\uncirc(\IDFT(\begin{bmatrix}
                     (C_1A_1B_1)^\dag &   &   \\
                       & \ddots   &   \\
                       &   &  (C_{n_3}A_{n_3}B_{n_3})^\dag \\
                   \end{bmatrix})))$.  $\Box$

\begin{theorem}
Let $\ac \in \mathbb{H}^{n_1\times n_2\times n_3}$,
$\dc \in \mathbb{H}^{l\times n_1\times n_3}$,
$\ec \in \mathbb{H}^{n_2\times k\times n_3}$.
If $\mathcal{A}$ is left invertible along $\mathcal{D},\mathcal{E}$, then
\begin{equation*}
  \mary{\mathcal{A}_l}{\mathcal{D},\mathcal{E}}=\mathcal{E}*(\mathcal{D}*\mathcal{A}*\mathcal{E})^\dag *\mathcal{D},
\end{equation*}
where $(\mathcal{D}*\mathcal{A}*\mathcal{E})^\dag=\fold(\uncirc(\IDFT(\begin{bmatrix}
                     (D_1A_1E_1)^\dag &   &   \\
                       & \ddots   &   \\
                       &   &  (D_{n_3}A_{n_3}E_{n_3})^\dag \\
                   \end{bmatrix})))$.
\end{theorem}
{\sc Proof:} Similar as Theorem \ref{tt5}. $\Box$


In the following, $A\in{H^{n_1\times n_2}_r}$ means that $A$ is a quaternion matrix with $rank(A)=r$. We will give some representations of the right/left inverse along two quaternion matrices/tensors by using the full rank decomposition.

\begin{theorem}\label{FG}{\rm\cite{Wu}} Let $A\in{H_r^{n_1\times n_2}}$. Then there exist
 matrices $F\in{H^{n_1\times r}}$ and $G\in{H^{r\times n_2}}$ with $rank(F)=rank(G)=r$
such that
\begin{equation}\label{PAPz}
  A=FG,
\end{equation}
\end{theorem}
which is called the full rank decomposition of $A$.

\begin{theorem}\label{tt7} Let $A\in{H^{n_1\times n_2}}$, $rank(A)=r$, $B \in {H}^{n_2\times l}$,
$C \in {H}^{k\times n_1}$. Suppose $B=\widehat{F}\widehat{G}$, $C=\widetilde{F}\widetilde{G}$, where $\widehat{F}\in{H_r^{n_2\times r}}$ and $\widehat{G}\in{H_r^{r\times l}}$, $\widetilde{F}\in{H_r^{k\times r}}$ and $\widetilde{G}\in{H_r^{r\times n_1}}$ are the full rank decomposition of $B$ and $C$, respectively. If $A$ is right invertible along $B,C$, then
\begin{equation*}
  \mary{A_r}{B,C}=\widehat{F}(\widetilde{G}A\widehat{F})^{-1}\widetilde{G}.
\end{equation*}
\end{theorem}
{\sc Proof:} Let us verify $X=\widehat{F}(\widetilde{G}A\widehat{F})^{-1}\widetilde{G}$ is the solution of (\ref{ATS2}). Firstly, it is easy to see
\begin{equation*}
\widehat{F}(\widetilde{G}A\widehat{F})^{-1}\widetilde{G}A\widehat{F}(\widetilde{G}A\widehat{F})^{-1}\widetilde{G}
=\widehat{F}(\widetilde{G}A\widehat{F})^{-1}\widetilde{G},
\end{equation*}
that is $\widehat{F}(\widetilde{G}A\widehat{F})^{-1}\widetilde{G}$ a $\{2\}$-inverse of $A$. Notice that
\begin{equation*}
R_r(\widehat{F}(\widetilde{G}A\widehat{F})^{-1}\widetilde{G})=R_r(X)\subseteq R_r(\widehat{F}), \qquad N_r(\widehat{F}(\widetilde{G}A\widehat{F})^{-1}\widetilde{G})=N_r(X)\supseteq N_r(\widetilde{G}).
\end{equation*}
Next, we will show the reverse part. As
\begin{equation*}
 \widehat{F}=\widehat{F}(\widetilde{G}A\widehat{F})^{-1}\widetilde{G}A\widehat{F}=
 \widehat{F}(\widetilde{G}A\widehat{F})^{-1}\widetilde{G}(A\widehat{F}),
\end{equation*}
we have $R_r(\widehat{F})\subseteq R(\widehat{F}(\widetilde{G}A\widehat{F})^{-1}\widetilde{G})=R_r(X)$.

Now, let $x\in N_r(\widehat{F}(\widetilde{G}A\widehat{F})^{-1}\widetilde{G})$, which implies $\widehat{F}(\widetilde{G}A\widehat{F})^{-1}\widetilde{G}x=0$. Pre-multiplying the equality by $(\widetilde{G}A\widehat{F})^{-1}\widetilde{G}A$, one has $(\widetilde{G}A\widehat{F})^{-1}\widetilde{G}A\widehat{F}(\widetilde{G}A\widehat{F})^{-1}\widetilde{G}x
=(\widetilde{G}A\widehat{F})^{-1}\widetilde{G}x=0$.  Again, pre-multiplication with $\widetilde{G}A\widehat{F}$ lead to $\widetilde{G}A\widehat{F}(\widetilde{G}A\widehat{F})^{-1} \widetilde{G}x=\widetilde{G}x=0$, which means $ N_r(\widehat{F}(\widetilde{G}A\widehat{F})^{-1}\widetilde{G})=N_r(X)\subseteq N_r(\widetilde{G})$. Therefore, $\mary{A_r}{B,C}=\widehat{F}(\widetilde{G}A\widehat{F})^{-1}\widetilde{G}.$ $\Box$

\begin{theorem} Let $A\in{H^{n_1\times n_2}}$, $rank(A)=r$, $D \in {H}^{l\times n_1}$,
$E \in {H}^{n_2\times k}$. Suppose $D=\widehat{F}\widehat{G}$, $E=\widetilde{F}\widetilde{G}$, where $\widehat{F}\in{H_r^{l\times r}}$ and $\widehat{G}\in{H_r^{r\times n_1}}$, $\widetilde{F}\in{H_r^{n_2\times r}}$ and $\widetilde{G}\in{H_r^{r\times k}}$ are the full rank decomposition of $D$ and $E$, respectively. If $A$ is left invertible along $D,E$, then
\begin{equation*}
  \mary{A_l}{D,E}=\widetilde{G}(\widehat{F}A\widetilde{G})^{-1}\widehat{F}.
\end{equation*}
\end{theorem}
{\sc Proof:} Similar as Theorem \ref{tt7}. $\Box$

By Theorem \ref{FG}, we can get the following definition. Notice that not all tensors have the full rank decomposition.
\begin{definition}\label{FG2} Let $\mathcal{A}\in \mathbb{H}^{n_1\times n_2\times n_3}$. If $\mathcal{A}$ can be decomposed into
$$\mathcal{A}=\mathcal{F}*\mathcal{G}$$
where $$\mathcal{F}=\fold({\uncirc}(\IDFT(\begin{bmatrix}
                     F_1 &   &   \\
                       & \ddots   &   \\
                       &   &  F_{n_3} \\
                   \end{bmatrix})))\in\mathbb{H}^{n_1\times r\times n_3},$$
and$$\mathcal{G}=\fold({\uncirc}(\IDFT(\begin{bmatrix}
                     G_1 &   &   \\
                       & \ddots   &   \\
                       &   &  G_{n_3} \\
                   \end{bmatrix})))\in\mathbb{H}^{r\times n_2 \times n_3},$$ $F_i\in{H^{n_1\times r}_{r}}$, $G_i\in{H^{r\times n_2}_{r}}$, $i=1,...,{n_3}$,
then we call the decomposition the full rank decomposition of $\mathcal{A}$.
\end{definition}

\begin{theorem}\label{tt8}
Let $\ac \in \mathbb{H}^{n_1\times n_2\times n_3}$,
$\bc \in \mathbb{H}^{n_2\times l\times n_3}$,
$\cc \in \mathbb{H}^{k\times n_1\times n_3}$.
Suppose $\mathcal{B}=\widehat{\mathcal{F}}*\widehat{\mathcal{G}}$, $\mathcal{C}=\widetilde{\mathcal{F}}*\widetilde{\mathcal{G}}$ are the full rank decomposition of $\mathcal{B}$ and $\mathcal{C}$, respectively. If $\mathcal{A}$ is right invertible along $\mathcal{B},\mathcal{C}$, then
\begin{equation*}
  \mary{\mathcal{A}_r}{\mathcal{B},\mathcal{C}}=\widehat{\mathcal{F}}*
  (\widetilde{\mathcal{G}}*\mathcal{A}*\widehat{\mathcal{F}})^{-1}*\widetilde{\mathcal{G}},
\end{equation*}
where $(\widetilde{\mathcal{G}}*\mathcal{A}*\widehat{\mathcal{F}})^{-1}=\fold(\uncirc(\IDFT(\begin{bmatrix}
                     (\widetilde{G}_1A_1\widehat{F}_1)^{-1} &   &   \\
                       & \ddots   &   \\
                       &   &  (\widetilde{G}_{n_3}A_{n_3}\widehat{F}_{n_3})^{-1} \\
                   \end{bmatrix})))$.
\end{theorem}
{\sc Proof:} Suppose \begin{eqnarray*} 
\DFT(\ccc(\unfold(\mathcal{A}))) =
\begin{bmatrix} A_1 & & \\
 & \ddots & \\
 & & A_{n_3}\end{bmatrix}.
\end{eqnarray*}
According to Theorem \ref{tt7}, we have
\begin{eqnarray*}
\begin{bmatrix} \mary{{A_1}_r}{B,C} & & \\
 & \ddots & \\
 & & \mary{{A_{n_3}}_r}{B_1,C_1} \end{bmatrix}&=&\begin{bmatrix} \widehat{F}_1(\widetilde{G}_1A_1\widehat{F}_1)^{-1} \widetilde{G}_1& & \\
 & \ddots & \\
 & & \widehat{F}_{n_3}(\widetilde{G}_{n_3}A_{n_3}\widehat{F}_{n_3})^{-1} \widetilde{G}_{n_3} \end{bmatrix}\\
 &=&\begin{bmatrix}
                      \widehat{F}_1 &   &   \\
                       & \ddots   &   \\
                       &   &  \widehat{F}_{n_3} \\
                   \end{bmatrix}\begin{bmatrix}
                     (\widetilde{G}_1A_1\widehat{F}_1)^{-1} &   &   \\
                       & \ddots   &   \\
                       &   &  (\widetilde{G}_{n_3}A_{n_3}\widehat{F}_{n_3})^{-1} \\
                   \end{bmatrix}\begin{bmatrix}
                     \widetilde{G}_1 &   &   \\
                       & \ddots   &   \\
                       &   &  \widetilde{G}_{n_3}\\
                   \end{bmatrix}.
\end{eqnarray*}
Performing $\fold(\text{uncirc}(\textbf{IDFT}))(\cdot)$ on the quaternion tensors of the both sides of the equation, one has
\begin{equation*}
  \mary{\mathcal{A}_r}{\mathcal{B},\mathcal{C}}=\widehat{\mathcal{F}}*
  (\widetilde{\mathcal{G}}*\mathcal{A}*\widehat{\mathcal{F}})^{-1}*\widetilde{\mathcal{G}},
\end{equation*}
where $(\widetilde{\mathcal{G}}*\mathcal{A}*\widehat{\mathcal{F}})^{-1}=\fold(\uncirc(\IDFT(\begin{bmatrix}
                     (\widetilde{G}_1A_1\widehat{F}_1)^{-1} &   &   \\
                       & \ddots   &   \\
                       &   &  (\widetilde{G}_{n_3}A_{n_3}\widehat{F}_{n_3})^{-1} \\
                   \end{bmatrix})))$. $\Box$

\begin{theorem}\label{tt9}
Let $\ac \in \mathbb{H}^{n_1\times n_2\times n_3}$,
$\dc \in \mathbb{H}^{l\times n_1\times n_3}$,
$\ec \in \mathbb{H}^{n_2\times k\times n_3}$.
Suppose $\mathcal{D}=\widehat{\mathcal{F}}*\widehat{\mathcal{G}}$, $\mathcal{E}=\widetilde{\mathcal{F}}*\widetilde{\mathcal{G}}$ are the full rank decomposition of $\mathcal{D}$ and $\mathcal{E}$, respectively. If $\mathcal{A}$ is left invertible along $\mathcal{D},\mathcal{E}$, then
\begin{equation*}
  \mary{\mathcal{A}_l}{\mathcal{D},\mathcal{E}}=\widetilde{\mathcal{G}}*
  (\widehat{\mathcal{F}}*\mathcal{A}*\widetilde{\mathcal{G}})^{-1}*\widehat{\mathcal{F}},
\end{equation*}
where $(\widehat{\mathcal{F}}*\mathcal{A}*\widetilde{\mathcal{G}})^{-1}=\fold(\uncirc(\IDFT(\begin{bmatrix}
                     (\widehat{F}_1A_1\widetilde{G}_1)^{-1} &   &   \\
                       & \ddots   &   \\
                       &   &  (\widehat{F}_{n_3}A_{n_3}\widetilde{G}_{n_3})^{-1} \\
                   \end{bmatrix})))$.
\end{theorem}
{\sc Proof:} Similar as Theorem \ref{tt8}. $\Box$

\section{An Algorithm for Computing the Inverse Along Two Quaternion Tensors}

In this part, we will give an algorithm for computing the right inverse along two quaternion tensors. Similarly, we can get the algorithm for computing the left inverse along two quaternion tensors. Examples are also given to demonstrate the algorithm.

\begin{algorithm}[H]
\caption{\textsc{Compute the right inverse of $\mathcal{A}$ along two quaternion tensors }}\label{algo:sylvester2}
\KwIn{$n_1\times n_2\times n_3$ quaternion tensor $\mathcal{A}$, $n_2\times l\times n_3$ quaternion tensor $\mathcal{B}$, $ k\times n_1\times n_3$ quaternion tensor $\cc$}

\KwOut{$n_2\times n_1\times n_3$ quaternion tensor $\mathcal{X}$}
\begin{enumerate}
\addtolength{\itemsep}{-0.8\parsep minus 0.8\parsep}

\item $i=3$, $\mathcal{D} = \fft(\mathcal{A},[ ~ ],i)$;

\item $i=3$, $\mathcal{E} = \fft(\mathcal{B},[ ~ ],i)$;

\item $i=3$, $\mathcal{F} = \fft(\mathcal{C},[ ~ ],i)$;

\item for $i=1,\ldots, n_3$

\quad $\mathcal{G}(:,:,i)=\mathcal{F}(:,:,i)\ast\mathcal{D}(:,:,i)\ast\mathcal{E}(:,:,i)$;

\quad $\mathcal{H}(:,:,i)= {\rm pinv}(\mathcal{G}(:,:,i))$, where
${\rm pinv}(\mathcal{G}(:,:,i))$ is the Moore-Penrose inverse of $\mathcal{G}(:,:,i)$;

\quad $\mathcal{T}(:,:,i)=\mathcal{E}(:,:,i)\ast\mathcal{H}(:,:,i)\ast\mathcal{F}(:,:,i)$;

end

\item  $i=3$, $\mathcal{X}=\ifft(\mathcal{T},[ ~ ],i)$.
\end{enumerate}
\end{algorithm}

\begin{myexample}
Let $\mathcal{A}$ be a $3\times3\times4$ quaternion tensor with the following form:

\begin{equation*}
  \mathcal{A}(:,:,1)=\begin{bmatrix}
                             1+{\bf i} &{\bf j} & 3\\
                             {\bf k} & 1+{\bf j} & {\bf i}+{\bf j}\\
                             2  & {\bf j}-{\bf k} & 1 \\
                          \end{bmatrix}, \qquad \mathcal{A}(:,:,2)=\begin{bmatrix}
                              2{\bf k} &{\bf i}+{\bf k} & -{\bf k}\\
                             {\bf i}+{\bf j}+2{\bf k} & 2+{\bf i} & {\bf i}-{\bf j}\\
                             3  & 2{\bf k}  & {\bf i}+{\bf k} \\
                          \end{bmatrix},
\end{equation*}

\begin{equation*}
  \mathcal{A}(:,:,3)=\begin{bmatrix}
                              {\bf j} &{\bf i}+{\bf k} & 2+{\bf i}\\
                             1+{\bf k} & {\bf j} & 2-{\bf j}\\
                             -{\bf i}  & 1-{\bf k} & 4+{\bf i}-{\bf k} \\
                          \end{bmatrix}, \qquad  \mathcal{A}(:,:,4)=\begin{bmatrix}
                              2{\bf i}+{\bf k} & 1 & 2+3{\bf i}+{\bf j}+2{\bf k}\\
                             5{\bf i}+2{\bf j} & -3-{\bf i}-2{\bf k} & 1-{\bf i}\\
                             2{\bf j}  & -{\bf i}+ 2{\bf j}+{\bf k} & 3+2{\bf k} \\
                          \end{bmatrix}.
\end{equation*}
Then, we apply Algorithm \ref{algo:sylvester2} to compute the Drazin inverse  of $\mathcal{A}$. In this case, $\mathcal{B}$ and $\mathcal{C}$ defined in Theorem \ref{tt5} hold $\mathcal{B}=\mathcal{C}=\mathcal{A}^k$. Hence, we have
\footnotesize{
\begin{align}
\nonumber
\begin{array}{lc}
\ac^D(:,:,1)  = \vspace{0.5ex}\\
    \begin{bmatrix}
   0.0342+0.0513{\bf i}-0.0940{\bf j}-0.0256{\bf k}  &   0.0726-0.1496{\bf i}+0.0385{\bf j}-0.0470{\bf k} &   -0.1581+0.0214{\bf i}+0.0299{\bf j}-0.0556{\bf k} \\
    0.0470+0.0983{\bf i}+0.0556{\bf j}-0.0556{\bf k}  &   -0.0128-0.1581{\bf i}-0.1838{\bf j}+0.0299{\bf k} &   -0.0214-0.0043{\bf i}-0.0385{\bf j}-0.1410{\bf k} \\
   -0.0556+0.0556{\bf i}+0.0385{\bf j}+0.0812{\bf k}  &   -0.0214+0.0043{\bf i}-0.0812{\bf j}-0.1838{\bf k} &   -0.0897+0.1068{\bf i}+0.2350{\bf j}+0.0641{\bf k} \\
                          \end{bmatrix},
\end{array}
\end{align}}
\footnotesize{
\begin{align}
\nonumber
\begin{array}{lc}
\ac^D(:,:,2)  = \vspace{0.5ex}\\
    \begin{bmatrix}
   0.0813-0.0255{\bf i}+0.0222{\bf j}-0.0539{\bf k}  &   0.0293+0.0927{\bf i}+0.0317{\bf j}+0.0241{\bf k} &   0.0284+0.0804{\bf i}-0.0099{\bf j}+0.0270{\bf k} \\
    0.0137-0.0407{\bf i}+0.0038{\bf j}-0.0345{\bf k}  &   0.0433-0.1411{\bf i}-0.0707{\bf j}-0.0404{\bf k} &   0.0714-0.0842{\bf i}+0.0719{\bf j}+0.0288{\bf k} \\
   -0.0440+0.0293{\bf i}+0.0066{\bf j}-0.0005{\bf k}  &   -0.0463+0.0648{\bf i}-0.0596{\bf j}-0.0468{\bf k} &   -0.1132+0.0158{\bf i}-0.0173{\bf j}-0.0494{\bf k} \\
                          \end{bmatrix},
\end{array}
\end{align}}
\footnotesize{
\begin{align}
\nonumber
\begin{array}{lc}
\ac^D(:,:,3)  = \vspace{0.5ex}\\
    \begin{bmatrix}
   0.0276-0.0104{\bf i}-0.0136{\bf j}-0.0548{\bf k}  &   0.0485+0.0318{\bf i}-0.0506{\bf j}-0.0929{\bf k} &   -0.1028-0.0167{\bf i}-0.0423{\bf j}-0.0276{\bf k} \\
    0.0211+0.0143{\bf i}-0.0206{\bf j}-0.0420{\bf k}  &   0.0451+0.0655{\bf i}+0.0029{\bf j}+0.1009{\bf k} &   -0.0271-0.0438{\bf i}-0.0151{\bf j}+0.0037{\bf k} \\
   -0.0550+0.0462{\bf i}+0.0013{\bf j}-0.0284{\bf k}  &   -0.0412+0.1069{\bf i}-0.1080{\bf j}-0.0777{\bf k} &   -0.1135+0.0222{\bf i}-0.0592{\bf j}+0.0827{\bf k} \\
                          \end{bmatrix},
\end{array}
\end{align}}
\footnotesize{
\begin{align}
\nonumber
\begin{array}{lc}
\ac^D(:,:,4)  = \vspace{0.5ex}\\
    \begin{bmatrix}
   0.0132-0.0264{\bf i}-0.0184{\bf j}-0.0084{\bf k}  &   0.0159-0.0167{\bf i}+0.0125{\bf j}-0.0194{\bf k} &   0.0437-0.0159{\bf i}+0.0138{\bf j}+0.0234{\bf k} \\
    -0.0030+0.0024{\bf i}-0.0092{\bf j}-0.0056{\bf k}  &   -0.0191+0.0059{\bf i}+0.0238{\bf j}-0.0097{\bf k} &   -0.0657+0.0166{\bf i}+0.0205{\bf j}+0.0330{\bf k} \\
   0.0135-0.0270{\bf i}+0.0110{\bf j}+0.0417{\bf k}  &   0.0548+0.0189{\bf i}+0.0183{\bf j}-0.0231{\bf k} &   -0.0080-0.0401{\bf i}+0.0298{\bf j}-0.0429{\bf k} \\
                          \end{bmatrix}.
\end{array}
\end{align}}
\end{myexample}

\begin{myexample}
Let $\ac \in \mathbb{H}^{3\times3\times4}$, $\bc \in \mathbb {H}^{3\times2\times4}$,
$\cc \in \mathbb {H}^{3\times3\times4}$ have the following forms. It is easy to check $\mathcal{A}$ is right invertible along $\mathcal{B},\mathcal{C}$.

\begin{equation*}
  \mathcal{A}(:,:,1)=\begin{bmatrix}
                             1+{\bf k} &{\bf i}+2{\bf j} & -{\bf i}-2{\bf k}\\
                             2 & 2+{\bf i}-{\bf k} & -2{\bf i}-{\bf j}\\
                             {\bf i}+{\bf j}  & 2{\bf k} & 1 \\
                          \end{bmatrix}, \qquad \mathcal{A}(:,:,2)=\begin{bmatrix}
                              3-{\bf i}-{\bf k} &{\bf i} & 1+{\bf i}\\
                             3{\bf j}+2{\bf k} & -3{\bf j}-{\bf k} & 1\\
                             {\bf j}  & 2+{\bf i}  & 2{\bf i}+2{\bf j} \\
                          \end{bmatrix},
\end{equation*}
\begin{equation*}
  \mathcal{A}(:,:,3)=\begin{bmatrix}
                              {\bf i} &{\bf j} & {\bf i}-{\bf j}\\
                             {\bf j}+3{\bf k} & 2 & 1+3{\bf i}\\
                             2{\bf j}-{\bf k}  & 1+{\bf i} & {\bf j} \\
                          \end{bmatrix}, \qquad  \mathcal{A}(:,:,4)=\begin{bmatrix}
                              {\bf i}+{\bf j} & {\bf k} & 2+{\bf i} \\
                             3{\bf i}  & 1  & 2+2{\bf k}\\
                             2{\bf i}-2{\bf j}  &  2{\bf j}-{\bf k} & 2{\bf j} \\
                          \end{bmatrix}.
\end{equation*}

\begin{equation*}
  \mathcal{B}(:,:,1)=\begin{bmatrix}
                             1   & {\bf i}+{\bf j} \\
                             {\bf k} & 1-{\bf j} \\
                             {\bf i}+{\bf j}+2{\bf k}  & 2-{\bf i}  \\
                          \end{bmatrix}, \qquad \mathcal{B}(:,:,2)=\begin{bmatrix}
                              2{\bf i}+{\bf j} & -{\bf i}+{\bf j}-2{\bf k} \\
                               2  & {\bf j}+{\bf k} \\
                             {\bf i}  & -2{\bf j}   \\
                          \end{bmatrix},
\end{equation*}
\begin{equation*}
  \mathcal{B}(:,:,3)=\begin{bmatrix}
                              {\bf k} &{\bf i}+{\bf k} \\
                             1-2{\bf j}  & 2{\bf i}+{\bf j} \\
                             1-{\bf i}-{\bf k}  & {\bf j}  \\
                          \end{bmatrix}, \qquad  \mathcal{B}(:,:,4)=\begin{bmatrix}
                              1+{\bf i}+{\bf j}-2{\bf k} & 2{\bf i}  \\
                             {\bf j}-{\bf k}  & 3  \\
                             2+{\bf k}  &  {\bf k}  \\
                          \end{bmatrix}.
\end{equation*}

\begin{equation*}
  \mathcal{C}(:,:,1)=\begin{bmatrix}
                             1  &{\bf i}+{\bf j} & {\bf k}\\
                             {\bf j}-{\bf k} & 2+{\bf i} & {\bf i}-{\bf j}\\
                             3-{\bf j}  & {\bf j}+{\bf k} & {\bf i} \\
                          \end{bmatrix}, \qquad \mathcal{C}(:,:,2)=\begin{bmatrix}
                              2+{\bf i}-{\bf j} &{\bf i} & {\bf j}\\
                             2{\bf i}-{\bf j} & 3+{\bf k} & 3{\bf i}\\
                             2{\bf j}+{\bf k}  & 1  & 2{\bf i}-{\bf k} \\
                          \end{bmatrix},
\end{equation*}
\begin{equation*}
  \mathcal{C}(:,:,3)=\begin{bmatrix}
                              {\bf i} &{\bf j}-{\bf k} & 3+2{\bf k}\\
                             1+{\bf j} & 5 & j+2{\bf k}\\
                             {\bf j}  & {\bf i}+{\bf j}+{\bf k} & 1 \\
                          \end{bmatrix}, \qquad  \mathcal{C}(:,:,4)=\begin{bmatrix}
                              {\bf k} & 1+{\bf i} & {\bf j}+{\bf k} \\
                             2-{\bf i}-2{\bf j}+3{\bf k}  & 2{\bf j}-{\bf k}  & {\bf i}+{\bf k}\\
                             3+{\bf i}-2{\bf j}  &  {\bf i} & {\bf j} \\
                          \end{bmatrix}.
\end{equation*}
Now, we will compute the right inverse of $\mathcal{A}$ along $\mathcal{B},\mathcal{C}$. Implement Algorithm \ref{algo:sylvester2} on $\mathcal{A}$, we have
\footnotesize{
\begin{align}
\nonumber
\begin{array}{lc}
\mary{\mathcal{A}_r}{\mathcal{B},\mathcal{C}}(:,:,1)  = \vspace{0.5ex}\\
    \begin{bmatrix}
   0.01278+0.0149{\bf i}-0.0156{\bf j}+0.1356{\bf k}  &   0.0588+0.0279{\bf i}-0.0269{\bf j}-0.0939{\bf k} &   0.0237+0.0310{\bf i}-0.0553{\bf j}+0.0163{\bf k} \\
    0.1406-0.0375{\bf i}-0.0301{\bf j}-0.0296{\bf k}  &   -0.0918+0.0275{\bf i}-0.1121{\bf j}+0.0527{\bf k} &   0.1054-0.0506{\bf i}-0.0246{\bf j}+0.0004{\bf k} \\
   -0.0730-0.0388{\bf i}-0.0003{\bf j}-0.0858{\bf k}  &   -0.0004+0.0148{\bf i}-0.0493{\bf j}+0.1090{\bf k} &   0.0226+0.0194{\bf i}+0.0202{\bf j}+0.0082{\bf k} \\
                          \end{bmatrix},
\end{array}
\end{align}}
\footnotesize{
\begin{align}
\nonumber
\begin{array}{lc}
\mary{\mathcal{A}_r}{\mathcal{B},\mathcal{C}}(:,:,2)  = \vspace{0.5ex}\\
    \begin{bmatrix}
   -0.0270-0.0271{\bf i}-0.0798{\bf j}-0.1330{\bf k}  &   0.1046-0.0362{\bf i}-0.0252{\bf j}+0.0646{\bf k} &   0.0486+0.0009{\bf i}+0.0381{\bf j}+0.0431{\bf k} \\
    -0.0364-0.0271{\bf i}-0.0900{\bf j}+0.0369{\bf k}  &   0.1031-0.0490{\bf i}+0.0275{\bf j}+0.0072{\bf k} &   0.0660-0.0221{\bf i}-0.0491{\bf j}-0.0041{\bf k} \\
   0.0314+0.0998{\bf i}-0.0239{\bf j}+0.0646{\bf k}  &   -0.0193-0.0337{\bf i}-0.0987{\bf j}-0.0014{\bf k} &   -0.0105+0.0540{\bf i}-0.0165{\bf j}-0.0092{\bf k} \\
                          \end{bmatrix},
\end{array}
\end{align}}
\footnotesize{
\begin{align}
\nonumber
\begin{array}{lc}
\mary{\mathcal{A}_r}{\mathcal{B},\mathcal{C}}(:,:,3)  = \vspace{0.5ex}\\
    \begin{bmatrix}
   0.1102-0.0758{\bf i}+0.0706{\bf j}-0.0288{\bf k}  &   0.0259+0.0098{\bf i}-0.0635{\bf j}+0.0477{\bf k} &   0.1233-0.0249{\bf i}+0.0981{\bf j}+0.0255{\bf k} \\
    -0.1346-0.0269{\bf i}+0.0370{\bf j}+0.1605{\bf k}  &   0.1292-0.0017{\bf i}-0.0685{\bf j}-0.0711{\bf k} &   0.0159-0.0380{\bf i}+0.0123{\bf j}-0.0962{\bf k} \\
   -0.0269-0.0768{\bf i}-0.0987{\bf j}+0.0685{\bf k}  &   0.0318+0.0237{\bf i}-0.0124{\bf j}-0.0181{\bf k} &   0.0274-0.1399{\bf i}-0.0648{\bf j}+0.0646{\bf k} \\
                          \end{bmatrix},
\end{array}
\end{align}}
\footnotesize{
\begin{align}
\nonumber
\begin{array}{lc}
\mary{\mathcal{A}_r}{\mathcal{B},\mathcal{C}}(:,:,4)  = \vspace{0.5ex}\\
    \begin{bmatrix}
   0.0528-0.0211{\bf i}-0.0203{\bf j}+0.0797{\bf k}  &   0.0774+0.0347{\bf i}+0.0206{\bf j}-0.1083{\bf k} &   0.0606-0.0479{\bf i}+0.1029{\bf j}+0.0408{\bf k} \\
    0.0763-0.1252{\bf i}-0.0851{\bf j}-0.1742{\bf k}  &   -0.0769-0.0066{\bf i}-0.0789{\bf j}+0.0786{\bf k} &   -0.0062-0.1280{\bf i}+0.0197{\bf j}+0.0126{\bf k} \\
   -0.0318-0.0723{\bf i}-0.1126{\bf j}-0.0060{\bf k}  &   0.0250-0.0268{\bf i}-0.0392{\bf j}+0.0876{\bf k} &   0.0167-0.1183{\bf i}-0.0480{\bf j}+0.0140{\bf k} \\
                          \end{bmatrix}.
\end{array}
\end{align}}
\end{myexample}

~

{\bf\large Funding}

~

This work was supported by Talent Introduction and Scientific Research Start-Up Project of Guangxi Minzu University (No.2021KJQD02); Guangxi Science and Technology Base and Talents Special Project (No.GUIKE21220024); Guangxi Natural Science Foundation (No.2018GXNSFDA281023) and National Natural Science Foundation of China (No.12061015).

\end{document}